\documentclass[a4paper, final, 12pt]{article}
\usepackage{amsmath, amssymb, latexsym, amscd, amsthm,amsfonts,amstext}
\usepackage[mathscr]{eucal}
\usepackage{graphicx}
\usepackage{subfig}
\usepackage{float}
\usepackage{color}
\usepackage{hyperref}
\usepackage[utf8]{inputenc}
\usepackage[english]{babel}
\usepackage[margin=1in]{geometry}

\numberwithin{equation}{section}
\input{tcilatex}
\begin{document}

\title{An Inverse Problem With the Final Overdetermination for the Mean
Field Games System}
\author{Michael V. Klibanov \thanks{
Department of Mathematics and Statistics, University of North Carolina at
Charlotte, Charlotte, NC, 28223, USA, mklibanv@uncc.edu} \and Jingzhi Li 
\thanks{
Department of Mathematics \& National Center for Applied Mathematics
Shenzhen \& SUSTech International Center for Mathematics, Southern
University of Science and Technology, Shenzhen 518055, P.~R.~China,
li.jz@sustech.edu.cn}, \and Hongyu Liu \thanks{
Department of Mathematics, City University of Hong Kong, Kowloon, Hong Kong
SAR, P.R. China, hongyliu@cityu.edu.hk}}
\date{}
\maketitle

\begin{abstract}
The mean field games (MFG)\ theory has broad application in mathematical
modeling of social phenomena. The Mean Field Games System (MFGS) is the key
to the MFG theory. This is a system of two nonlinear parabolic partial
differential equations with two opposite directions of time $t\in \left(
0,T\right) .$ The topic of Coefficient Inverse Problem (CIPs) for the MFGS
is a newly emerging one. A CIP for the MFGS is studied. The input data are
Dirichlet and Neumann boundary conditions either on a part of the lateral
boundary (incomplete data) or on the whole lateral boundary (complete data).
In addition to the initial conditions at $\left\{ t=0\right\} ,$ terminal
conditions at $\left\{ t=T\right\} $ are given. The terminal conditions mean
the final overdetermination. The necessity of assigning all these input data
is explained. H\"{o}lder and Lipschitz stability estimates are obtained for
the cases of incomplete and complete data respectively. These estimates
imply uniqueness of the CIP.
\end{abstract}

\textbf{Key Words}: the mean field games system, new Carleman estimates, H%
\"{o}lder and Lipschitz stability estimates, uniqueness,

\textbf{2020 MSC codes}: 35R30

\section{Introduction}

\label{sec:1}

Social sciences play a significant role in the modern society. The Mean
Field Games (MFG) theory studies the collective behavior of large
populations of rational decision-makers (agents). This theory has a number
of applications in the mathematical modeling of social phenomena. Some
examples of these applications are, e.g. finance \cite{A,Trusov}, sociology 
\cite{Bauso}, fighting corruption \cite{KM,Kol}, cyber security \cite{Kol},
etc.

This theory was first introduced in 2006-2007 in seminal works of Lasry and
Lions \cite{LL1,LL2,LL3} and of Huang, Caines and Malham\'{e} \cite%
{Huang1,Huang2}. The mean field games system (MFGS) is the core of the MFG
theory. This is a system of two coupled nonlinear parabolic Partial
Differential Equations (PDEs) with two opposite directions of time. In the
first equation time is running downwards. This is the
Hamilton-Jacobi-Bellman equation (HJB). And in the second equation time is
running upwards. This is Fokker-Planck (FP) equation. Let $\Omega \subset 
\mathbb{R}^{n},n\geq 1$ be a bounded domain with its boundary $\partial
\Omega $ and let time $t\in \left( 0,T\right) .$ The state position of a
representative agent is $x\in \Omega .$ HJB equation governs the value
function $u\left( x,t\right) $ of each individual agent located at $x$ at
the moment of time $t.$ FP equation describes the evolution of the
distribution of agents $m\left( x,t\right) $ over time $t\in \left(
0,T\right) $.

Due to the applications of the MFG theory, it is important to study a
variety of mathematical topics of this theory. In the current paper, we
study a Coefficient Inverse Problem (CIP) for the MFGS. We consider the case
of the data resulting from a single measurement event. Previously a H\"{o}%
lder stability estimate was obtained in \cite{MFG6} for a CIP for the MFGS
with a single measurement data. However, the statement of the CIP in \cite%
{MFG6} is significantly different from the one of this paper, see subsection
3.1.

CIPs for the MFGS is a newly emerging topic. \ We are aware only about six
previous publications about such CIPs, and we list them in this paragraph.
Stability and uniqueness theorems for CIPs for the MFGS with single
measurement data were obtained in \cite{ImLY,MFG6}. Uniqueness theorems for
the case of infinitely many measurements were obtained in \cite%
{Liu1,Liu2,Ren}. We refer to \cite{Chow,Ding} for numerical studies of CIPs
for the MFGS.

Since the input data for CIPs are results of measurements, then they are
given with errors. Hence, we are concerned with obtaining H\"{o}lder and
Lipschitz stability estimates of the solution of our CIP with respect to the
error in the input data. These stability estimates immediately imply
uniqueness of our CIP.

In this paper we modify the framework, which was first proposed in \cite%
{BukhKlib}, where the apparatus of Carleman estimates was introduced in the
field of CIPs, also, see, e.g. \cite%
{ImYam1,Isakov,Klib84,Klib92,Ksurvey,KL,Yam} and references cited therein
for some follow up publications. Carleman estimates were introduced in the
MFG theory in \cite{MFG1} and were used since then in \cite{ImLY,MFG4,MFG6}
as well as in the current paper. The idea of our modification of the
framework of \cite{BukhKlib} is outlined in subsection 3.2.

It is natural to call the problem of this paper \textquotedblleft CIP with
the final overdetermination". Indeed, we assume that we know both initial
and terminal conditions for both functions $u$ and $m$ as well as both\
Dirichlet and Neumann boundary conditions for these functions on either a
part of the lateral boundary or on the whole lateral boundary. On the other
hand, if a CIP for a single parabolic equation requires to find a
coefficient of this equation, assuming that its solution is known at $%
\left\{ t=0\right\} $ and at $\left\{ t=T\right\} ,$ then such a CIP is
called \textquotedblleft CIP\ with the final overdetermination": we refer to 
\cite[section 9.1]{Isakov} for the Lipschitz stability result for such a
problem for the case of a single parabolic equation. However, the case of
the MFGS with the final overdetermination was not considered previously. In
addition, the technique of this paper is significantly different from the
one of \cite{Isakov}. A more detailed discussion of the statement of our CIP
can be found in subsections 3.1 and 3.2.

\textbf{Remark 1.1}. \emph{Traditionally minimal smoothness assumptions are
of a secondary concern in the field of inverse problems, see, e.g. \cite{Nov}%
, \cite[Theorem 4.1]{Rom}. Therefore, we are not concerned with these
assumptions in the current paper. }

All functions below are real valued ones. In section 2, we first formulate
the MFGS and outline four main difficulties of working with this system.
Then we formulate our CIP. In section 3, we first discuss our input data for
our CIP. Next, we outline our idea of the above mentioned modification of
the framework of \cite{BukhKlib}. We formulate our theorems in section 4.
Two theorems of this section about H\"{o}lder and Lipschitz stability are
proven in section 5. On the other hand, two theorems of section 4 about
Carleman estimates are proven in Appendix.

\section{Statement of the Coefficient Inverse Problem}

\label{sec:2}

\subsection{Domains and spaces}

\label{sec:2.1}

First, we introduce some basic notations we use in this paper. For $x=\left(
x_{1},x_{2},...,x_{n}\right) \in \mathbb{R}^{n}$ denote $\overline{x}=\left(
x_{2},...,x_{n}\right) .$ To simplify the presentation, we consider only the
case when our domain of interest $\Omega \subset \mathbb{R}^{n}$ is a
rectangular prism. Let $A_{i}>0,i=1,...,n$ and $T>0$ be some numbers. Also,
let $\gamma \in \left( 0,2A_{1}\right) $ be another number. Everywhere below
domains and their boundaries are defined as: 
\begin{equation}
\left. 
\begin{array}{c}
\Omega =\left\{ x:-A_{i}<x_{i}<A_{i},i=1,...,n\right\} ,\text{ } \\ 
\Omega ^{\prime }=\left\{ \overline{x}:-A_{i}<x_{i}<A_{i},i=2,...,n\right\} ,
\\ 
\Gamma ^{-}=\left\{ x\in \partial \Omega :x_{1}=-A_{1}\right\} ,\Gamma
^{+}=\left\{ x\in \partial \Omega :x_{1}=A_{1}\right\} ,\Gamma _{T}^{\pm
}=\Gamma ^{\pm }\times \left( 0,T\right) , \\ 
\partial _{i}^{\pm }\Omega =\left\{ x\in \partial \Omega :x_{i}=\pm
A_{i}\right\} ,\text{ }\partial _{i}^{\pm }\Omega _{T}=\partial _{i}^{\pm
}\Omega \times \left( 0,T\right) ,\text{ }i=2,...,n, \\ 
Q_{T}=\Omega \times \left( 0,T\right) ,\text{ }S_{T}=\partial \Omega \times
\left( 0,T\right) , \\ 
\Omega _{\gamma }=\left\{ x\in \Omega :x_{1}\in \left( -A_{1}+\gamma
,A_{1}\right) \right\} ,Q_{\gamma T}=\Omega _{\gamma }\times \left(
0,T\right) .%
\end{array}%
\right.  \label{2.1}
\end{equation}%
We now introduce some spaces we will work with. Let $k\geq 1$ be an
integers. Denote%
\begin{equation}
C^{2k,k}\left( \overline{Q}_{T}\right) =\left\{ u:\left\Vert u\right\Vert
_{C^{2k,k}\left( \overline{Q}_{T}\right) }=\max_{\left\vert \alpha
\right\vert +2m\leq 2k}\left\Vert D_{x}^{\alpha }\partial
_{t}^{m}u\right\Vert _{C\left( \overline{Q}_{T}\right) }<\infty \right\} ,
\label{2.100}
\end{equation}%
\begin{equation}
H^{2k,k}\left( Q_{T}\right) =\left\{ u:\left\Vert u\right\Vert
_{H^{4,2}\left( Q_{T}\right) }^{2}=\dsum\limits_{\left\vert \alpha
\right\vert +2m\leq 2k}\left\Vert D_{x}^{\alpha }\partial
_{t}^{m}u\right\Vert _{L_{2}\left( Q_{T}\right) }^{2}<\infty \right\} ,
\label{2.101}
\end{equation}%
\begin{equation}
H^{2,1}\left( \partial _{i}^{\pm }\Omega _{T}\right) =\left\{ 
\begin{array}{c}
u:\left\Vert u\right\Vert _{H^{2,1}\left( \partial _{i}^{\pm }\Omega
_{T}\right) }^{2}=\dsum\limits_{j=1,j\neq i}^{n}\left\Vert
u_{x_{j}}\right\Vert _{L_{2}\left( \partial _{i}^{\pm }\Omega _{T}\right)
}^{2}+ \\ 
+\dsum\limits_{j,s=1,\left( j,s\right) \neq \left( i,i\right)
}^{n}\left\Vert u_{x_{j}x_{s}}\right\Vert _{L_{2}\left( \partial _{i}^{\pm
}\Omega _{T}\right) }^{2}+ \\ 
+\dsum\limits_{j=0}^{1}\left\Vert \partial _{t}^{j}u\right\Vert
_{L_{2}\left( \partial _{i}^{\pm }\Omega _{T}\right) }^{2}<\infty%
\end{array}%
\right\} ,  \label{2.102}
\end{equation}%
\begin{equation}
H^{2,1}\left( S_{T}\right) =\left\{ u:\left\Vert u\right\Vert
_{H^{2,1}\left( S_{T}\right) }^{2}=\dsum\limits_{i=1}^{n}\left\Vert
u\right\Vert _{H^{2,1}\left( \partial _{i}^{\pm }\Omega _{T}\right)
}^{2}<\infty \right\} ,  \label{2.103}
\end{equation}%
\begin{equation}
H^{1,0}\left( \partial _{i}^{\pm }\Omega _{T}\right) =\left\{ u:\left\Vert
u\right\Vert _{H^{1,0}\left( \partial _{i}^{\pm }\Omega _{T}\right)
}^{2}=\dsum\limits_{j=1,j\neq i}^{n}\left\Vert u_{x_{j}}\right\Vert
_{L_{2}\left( \partial _{i}^{\pm }\Omega _{T}\right) }^{2}+\left\Vert
u\right\Vert _{L_{2}\left( \partial _{i}^{\pm }\Omega _{T}\right)
}^{2}<\infty \right\} ,  \label{2.104}
\end{equation}%
\begin{equation}
H^{1,0}\left( S_{T}\right) =\left\{ u:\left\Vert u\right\Vert
_{H^{1,0}\left( S_{T}\right) }^{2}=\dsum\limits_{i=1}^{n}\left\Vert
u\right\Vert _{H^{1,0}\left( \partial _{i}^{\pm }\Omega _{T}\right)
}^{2}<\infty \right\} .  \label{2.105}
\end{equation}%
Spaces $H^{2,1}\left( \Gamma _{T}^{\pm }\right) $ and $H^{1,0}\left( \Gamma
_{T}^{\pm }\right) $ are defined similarly with spaces $H^{2,1}\left(
\partial _{i}^{\pm }\Omega _{T}\right) $ and $H^{1,0}\left( \partial
_{i}^{\pm }\Omega _{T}\right) $ respectively.

\subsection{The mean field games system}

\label{sec:2.2}

We consider a slightly simplified form of the MFGS of the second order  \cite%
{A,LL3}:%
\begin{equation}
\left. 
\begin{array}{c}
u_{t}(x,t)+\Delta u(x,t){-a(x)(\nabla u(x,t))^{2}/2}+ \\ 
+\dint\limits_{\Omega }K\left( x,y\right) m\left( y,t\right) dy+s\left(
x,t\right) m\left( x,t\right) =0,\text{ }\left( x,t\right) \in Q_{T}, \\ 
m_{t}(x,t)-\Delta m(x,t){-\func{div}(a(x)m(x,t)\nabla u(x,t))}=0,\text{ }%
\left( x,t\right) \in Q_{T},%
\end{array}%
\right.   \label{2.2}
\end{equation}%
where ${\nabla u=}\left( u_{x_{1}},...,u_{x_{n}}\right) ,$ and conditions on
functions $K\left( x,y\right) ,$ $a(x)$ and $s\left( x,t\right) $ are
imposed later.

The term with the integral in (\ref{2.2}) is called \textquotedblleft global
interaction term". This term has a deep applied meaning, which is explained
in \cite[page 634]{MFG1}. More precisely, $K\left( x,y\right) $ is the
action on the agent occupying the point $x$ by the agent occupying the point 
$y$. Hence, the integral term in (\ref{2.2}) expresses the average action of
all agents located at all points $y\in \Omega $ on the agent located at the
point $x$.

\subsection{Four main difficulties of working with MFGS (\protect\ref{2.2})}

\label{sec:2.3}

We now outline four main difficulties of working with MFGS (\ref{2.2}):

\begin{enumerate}
\item MFGS (\ref{2.2}) is highly nonlinear. On the top of this is that any
CIP for a PDE is nonlinear as well.

\item Two equations of system (\ref{2.2}) have two opposite directions of
time. Therefore, the classical theory of parabolic equations \cite{Ladpar}
does not work here.

\item The presence of the integral term in the first equation (\ref{2.2}).
Such terms are not present in all past works on CIPs for parabolic PDEs.

\item The presence of the Laplace operator $\Delta u$ in the second equation
(\ref{2.2}), since this operator is involved in the principal part $\left(
\partial _{t}+\Delta \right) u$ of the first equation (\ref{2.2}).
\end{enumerate}

Due to items 1-4, the past theory of CIPs for a single parabolic PDE cannot
be automatically applied to CIPs for MFGS (\ref{2.2}). Rather, a significant
additional effort is required for the latter.

\subsection{The kernel $K\left( x,y\right) $ in (\protect\ref{2.2})}

\label{sec:2.4}

Let $M>0$ be a number, $\delta \left( z\right) ,$ $z\in \mathbb{R}$ be the
delta-function and $H\left( z\right) $ be the Heaviside function, 
\begin{equation*}
H\left( z\right) =\left\{ 
\begin{array}{c}
1\text{ if }z>0, \\ 
0\text{ if }z<0.%
\end{array}%
\right. 
\end{equation*}%
We cannot work with a general function $K\left( x,y\right) .$ Hence, we
assume below that the integral term in the first equation (\ref{2.2}) has
the following form: 
\begin{equation}
\left. 
\begin{array}{c}
K\left( x,y\right) =b\left( x\right) \left\{ \delta \left(
x_{1}-y_{1}\right) K_{1}\left( \overline{x},\overline{y}\right) +H\left(
y_{1}-x_{1}\right) K_{2}\left( x,y\right) \right\} .\text{ } \\ 
K_{1}\in C^{4}\left( \overline{\Omega ^{\prime }}\times \overline{\Omega
^{\prime }}\right) ,\text{ }\left\Vert K_{1}\right\Vert _{C^{4}\left( 
\overline{\Omega ^{\prime }}\times \overline{\Omega ^{\prime }}\right) }<M,%
\text{ } \\ 
K_{2}\in C^{4}\left( \overline{\Omega }\times \overline{\Omega }\right) ,%
\text{ }\left\Vert K_{2}\right\Vert _{C^{4}\left( \overline{\Omega }\times 
\overline{\Omega }\right) }<M.%
\end{array}%
\right.   \label{2.3}
\end{equation}%
Hence, 
\begin{equation}
\left. 
\begin{array}{c}
\dint\limits_{\Omega }K\left( x,y\right) m\left( y,t\right) dy= \\ 
=b\left( x\right) \left[ \dint\limits_{\Omega ^{\prime }}K_{1}\left( 
\overline{x},\overline{y}\right) m\left( x_{1},\overline{y},t\right) d%
\overline{y}+\dint\limits_{x_{1}}^{A_{1}}\left( \dint\limits_{\Omega
^{\prime }}K_{2}\left( x,y_{1},\overline{y}\right) m\left( y_{1},\overline{y}%
,t\right) d\overline{y}\right) dy_{1}\right] .%
\end{array}%
\right.   \label{2.30}
\end{equation}%
Conditions imposed on the function $b\left( x\right) $ are specified later.
A popular example of $K\left( x,y\right) $ is \cite[section 4.2]{LiuOsher}: 
\begin{equation*}
K\left( x,y\right) =b\left( x\right) \frac{1}{\left( 2\pi \right) ^{n}}%
\dprod\limits_{i=1}^{n}\frac{1}{\sigma _{i}}\exp \left( -\frac{\left(
x_{i}-y_{i}\right) ^{2}}{2\sigma _{i}^{2}}\right) .
\end{equation*}%
We recall that Gaussians approximate the $\delta -$function in the sense of
distributions, which justifies our choice of $\delta \left(
x_{1}-y_{1}\right) K_{1}\left( \overline{x},\overline{y}\right) $ in (\ref%
{2.3}).\emph{\ }

\subsection{Coefficient Inverse Problem}

\label{sec:2.5}

It is hard to find a specific form of the kernel $K\left( x,y\right) ,$ see,
e.g. \cite[section 4]{LiuOsher}. Hence, the recovery of at least a part of
this kernel is of a significant interest, and this is what we do in the
current paper. More precisely, we are interested in this paper in the
recovery of the coefficient $b\left( x\right) $ in (\ref{2.3}). Following
Remark 1.1, we are not concerned here with some extra smoothness conditions
we impose below.

\textbf{Coefficient Inverse} \textbf{Problem }(CIP). \emph{Assume that
functions }$u,m\in C^{6,3}\left( \overline{Q}_{T}\right) $\emph{\ satisfy
equations (\ref{2.2}), and let condition (\ref{2.3}) holds. Let}%
\begin{equation}
\left. 
\begin{array}{c}
u\left( x,0\right) =p\left( x\right) ,\text{ }m\left( x,0\right) =q\left(
x\right) ,\text{ }x\in \Omega , \\ 
u\left( x,T\right) =F\left( x\right) ,\text{ }m\left( x,T\right) =G\left(
x\right) ,\text{ }x\in \Omega , \\ 
u\mid _{S_{T}}=f_{0}\left( x,t\right) ,\text{ }\partial _{n}u\mid
_{S_{T}}=f_{1}\left( x,t\right) , \\ 
m\mid _{S_{T}}=g_{0}\left( x,t\right) ,\text{ }\partial _{n}u\mid
_{S_{T}}=g_{1}\left( x,t\right) .%
\end{array}%
\right.  \label{2.4}
\end{equation}%
\emph{Determine the coefficient }$b\left( x\right) $\emph{, assuming that
the functions in the right hand sides of (\ref{2.4}) are known.}

Thus, the functions in the right hand sides of first two lines of (\ref{2.4}%
) are initial and terminal conditions. The right hand sides of the third and
fourth lines of (\ref{2.4}) are Dirichlet and Neumann boundary data, which
are also called \textquotedblleft lateral Cauchy data". We will consider the
following two cases of the lateral Cauchy data:

\begin{enumerate}
\item Incomplete lateral Cauchy data. This is the case when in (\ref{2.4}) 
\begin{equation}
\text{functions }f_{0},f_{1},g_{0},g_{1}\text{ are known at }S_{T}\diagdown
\Gamma _{T}^{-}\text{ and unknown at }\Gamma _{T}^{-}.  \label{2.5}
\end{equation}%
We obtain a H\"{o}lder stability estimate in this case.

\item Complete lateral Cauchy data. This is the case when in (\ref{2.4}) 
\begin{equation}
\text{functions }f_{0},f_{1},g_{0},g_{1}\text{ are known at the whole
boundary }S_{T}.  \label{2.6}
\end{equation}%
\emph{\ \ }We obtain Lipschitz stability estimate in this case.
\end{enumerate}

The input data (\ref{2.4}) are generated by a single measurement event. As
to these data, in the conventional case of the MFG theory, only functions $%
u\left( x,T\right) $ and $m\left( x,0\right) $ are given \cite{A} as well as
either Neumann or Dirichlet boundary condition for each of functions $u,m.$
In the case of a practical mean field game process, other functions in (\ref%
{2.4}) can be obtained via, e.g. polling of game participants at $t=0,T$ as
well as at the boundary $\partial \Omega ,$ see, e.g. \cite[page 2]{Chow}.

\section{Discussion}

\label{sec:3}

In this section we explain first why do we need the input data (\ref{2.4}).
Next, we briefly outline our idea of a modification of the framework of \cite%
{BukhKlib} in order to make it applicable to our CIP.

\subsection{Discussion of the input data (\protect\ref{2.4})}

\label{sec:3.1}

Let $t_{0}\in \left( 0,T\right) $ be a fixed moment of time. A H\"{o}lder
stability estimate was obtained in \cite{MFG6} for a CIP for MFGS (\ref{2.2}%
) in the case when the coefficient $a\left( x\right) $ is unknown, initial
and terminal conditions in (\ref{2.4}) 
\begin{equation}
u\left( x,0\right) ,u\left( x,T\right) ,m\left( x,0\right) ,m\left(
x,T\right)   \label{2.50}
\end{equation}%
are replaced with the assumption of the knowledge of functions $u\left(
x,t_{0}\right) $ and $m\left( x,t_{0}\right) ,$ and also lateral Cauchy data
in (\ref{2.4}) are known in \cite{MFG6}. In \cite{ImLY} Lipschitz stability
estimate was obtained for a similar CIP for MFGS (\ref{2.2}) without the
integral term in it.

In the case of a single parabolic equation, uniqueness and stability results
for CIPs with $x-$dependent unknown coefficients were obtained only under
the assumption that the solution of that equation is known at $t=t_{0}\in
\left( 0,T\right) $ and the lateral Cauchy data are known as well, see, e.g. 
\cite{ImYam1,Klib84,Klib92,Ksurvey,KL,Yam}. A similar statement is true for
CIPs for MFGS (\ref{2.2}) \cite{ImLY,MFG6}. These results were obtained
using the framework of \cite{BukhKlib}.

If, however, only the initial condition at $\left\{ t=0\right\} $ and
lateral Cauchy data are known, then a methodology of obtaining stability
results for such CIPs does not exist yet even for the case of a single
parabolic PDE. This explains our need of the knowledge in (\ref{2.4}) of all
four initial and terminal conditions (\ref{2.50}) as well as the lateral
Cauchy data.

On the other hand, if we would assume only the knowledge of functions (\ref%
{2.50}) as well as of only either Dirichlet or Neumann boundary condition
for each of functions $u,m$, then we would not be able to consider the case
of incomplete data (\ref{2.5}). In addition, we would likely need to impose
some yet unknown additional conditions on operators in (\ref{2.2}). For
example, a similar CIP for a single parabolic equation with the data at $%
\left\{ t=0\right\} $, $\left\{ t=T\right\} $ and the Dirichlet boundary
condition at the entire boundary is considereed in \cite[section 9.1]{Isakov}%
. And it is assumed in \cite{Isakov} that the Dirichlet boundary value
problem for the corresponding elliptic operator has no more than one
solution.

\subsection{Our modification of\ the\ framework\ of \protect\cite{BukhKlib}
\ }

\label{sec:3.2}

The first step the\ framework\ of \cite{BukhKlib} transforms the considered
CIP in an integral-differential equation, which does not contain the unknown
coefficient. Integral terms in this equation are $t-$dependent Volterra
integrals. Next, the application of a Carleman estimate to that equation
leads to the desired result. This scheme works in the case of CIPs for
hyperbolic and elliptic PDEs in the cases when the lateral Cauchy data are
given, in addition to some initial conditions, see, e.g. \cite%
{BY,Klib84,Klib92,Ksurvey,Yam}, \cite[Chapter 3]{KL}. And in the case of
CIPs for parabolic equations, this scheme works only if one replaces the
initial data at $\left\{ t=0\right\} $ with the data at $\left\{ t=t_{0}\in
\left( 0,T\right) \right\} ,$ see subsection 3.1.

However, this framework does not work for our case when both initial data at 
$\left\{ t=0\right\} $ and terminal data at $\left\{ t=T\right\} $ are given
in (\ref{2.4}). More precisely, the straightforward application of the
framework of \cite{BukhKlib} to our CIP leads to the presence of some
parasitic integrals over $\left\{ t=0\right\} $ and $\left\{ t=T\right\} .$
These integrals appear when integrating the pointwise Carleman estimate over
the time cylinder $Q_{T}$ and applying the Gauss formula. The presence of
these integrals does not allow us to obtain our desired H\"{o}lder and
Lipschitz stability estimates and uniqueness theorem for our CIP.

Hence, we modify here the idea of \cite{BukhKlib}. More precisely, we
arrange the above transformation in such a way that those parasitic
integrals cancel each other. After our transformation, each of two
transformed functions obtained from functions $u$ and $m$ in (\ref{2.2})
attains the same values at $\left\{ t=0\right\} $ and at $\left\{
t=T\right\} .$ Next, we apply to the resulting transformed system of
integral differential equations two new Carleman estimates for operators $%
\partial _{t}+\Delta $ and $\partial _{t}-\Delta $. The new point of these
estimates is that the Carleman Weight Function (CWF) in them is independent
on $t$. Our CWF depends only on $x_{1}$: due to (\ref{2.3}) and (\ref{2.30}%
). On the other hand, in conventional Carleman estimates for parabolic
equations with the lateral Cauchy data, CWFs always depend on both $x$ and $t
$, see, e.g. \cite[section 2.3]{KL}, \cite[\S 1 of Chapter 4]{LRS}, \cite%
{Yam}.

\section{Formulations of Theorems}

\label{sec:4}

\subsection{Carleman estimates}

\label{sec:4.1}

It is sufficient to prove Carleman estimates only for principal parts of
Partial Differential Operators \cite[Lemma 2.1.1]{KL}. There are two methods
of proofs of Carleman estimates. The first method is presented in books \cite%
[sections 8.3 and 8.4]{Horm}, \cite[Theorem 3.2.1]{Isakov}, and it is based
on symbols of Partial Differential Operators. This method is both elegant
and short. However, it is based on the assumption of zero boundary
conditions of involved functions. On the other hand, we work here with the
non-zero boundary conditions, which play an important role in stability
estimates for our CIP.

Therefore, primeraly due to our need to arrange the mutual cancellation of
parasitic integrals over $\left\{ t=0\right\} $ and $\left\{ t=T\right\} $
(see subsection 3.2), we need a painstaking analysis of boundary terms in
our Carleman estimates. Thus, we use the second method. By this method, one
first derives a pointwise Carleman estimate. Next, one integrates this
estimate over the domain of interest. Boundary integrals occur due to the
Gauss formula. In addition to our analysis of resulting integrals over $%
\left\{ t=0\right\} $ and $\left\{ t=T\right\} ,$ our derivation also allows
us to analyze resulting boundary integrals over the lateral boundary $S_{T},$
which is important for our target stability estimates of $H^{2,1}\left(
Q_{\gamma T}\right) $ and $H^{2,1}\left( Q_{T}\right) $ norms of involved
functions.

The derivation of any pointwise Carleman estimate is inevitably space
consuming, see, e.g. \cite[section 2.3]{KL}, \cite[\S 1 of Chapter 4]{LRS}
and \cite{Yam}. However, this is the price we pay for the incorporation of
non-zero boundary conditions.

We remind that due to (\ref{2.3}) and (\ref{2.30}), our CWF depends only on
the variable $x_{1}.$ On the other hand, as stated in subsection 3.2, it is
critical for our CIP that CWF should be independent on $t$, which is unusual
in Carleman estimates for parabolic operators.

Let $\nu >2$ and $\lambda >1$ be some large parameters, which we will choose
later. Consider two functions $\psi $ and $\varphi _{\lambda ,\nu },$ where $%
\varphi _{\lambda ,\nu }\left( x\right) $ is the CWF we work with. Thus, 
\begin{equation}
\psi \left( x\right) =x_{1}+A_{1}+2,\text{ }\varphi _{\lambda ,\nu }\left(
x\right) =e^{2\lambda \psi ^{\nu }},  \label{4.1}
\end{equation}%
\begin{equation}
\exp \left( 2\lambda \cdot 2^{\nu }\right) \leq \varphi _{\lambda ,\nu
}\left( x\right) \leq \exp \left[ 2\lambda \left( 2A_{1}+2\right) ^{\nu }%
\right] \text{in }\Omega .  \label{4.2}
\end{equation}

\textbf{Theorem 4.1 }(pointwise Carleman estimate for the operator $\partial
_{t}-\Delta $). \emph{There exist sufficiently large numbers }$\nu _{0}=\nu
_{0}\left( A_{1}\right) >2$\emph{, }$\lambda _{0}=\lambda _{0}\left(
A_{1}\right) >1$\emph{\ and a number }$C=C\left( Q_{T}\right) >0$ \emph{%
depending only on the domain }$\Omega $ \emph{such that \ for }$\nu =\nu
_{0},$ \emph{\ for all }$\lambda \geq \lambda _{0}$\emph{\ and for all
functions }$u\in C^{4,2}\left( \overline{Q}_{T}\right) $\emph{\ the
following pointwise Carleman estimate holds:\ }%
\begin{equation}
\left. 
\begin{array}{c}
\left( u_{t}-\Delta u\right) ^{2}\varphi _{\lambda ,\nu _{0}}\geq \left(
C/\lambda \right) \left(
u_{t}^{2}+\dsum\limits_{i,j=1}^{n}u_{x_{i}x_{j}}^{2}\right) \varphi
_{\lambda ,\nu _{0}}+ \\ 
+C\left[ \lambda \left( \nabla u\right) ^{2}+\lambda ^{3}u^{2}\right]
\varphi _{\lambda ,\nu _{0}}+\partial _{t}V+\func{div}U,\text{ }\left(
x,t\right) \in Q_{T},%
\end{array}%
\right.  \label{4.3}
\end{equation}%
\emph{where }$U$\emph{\ is an }$n-$\emph{D vector function. The function }$%
\partial _{t}V$\emph{\ is:}%
\begin{equation}
\left. 
\begin{array}{c}
\partial _{t}V= \\ 
=\partial _{t}\left[ \left( 2\lambda /\left( 2\lambda +1\right) \right)
\left( \left( u_{x_{1}}+\lambda \nu _{0}\psi ^{\nu _{0}-1}u\right)
^{2}+\dsum\limits_{i,j=2}^{n}u_{x_{i}}^{2}\right) \psi ^{-\nu _{0}+1}\varphi
_{\lambda ,\nu _{0}}\right] + \\ 
+\partial _{t}\left[ \left( 2\lambda /\left( 2\lambda +1\right) \right)
\left( -\lambda ^{2}\nu _{0}^{2}\psi ^{\nu _{0}-1}\left( 1-2\psi ^{-\nu
_{0}}\left( \nu _{0}-1\right) /\left( \lambda \nu _{0}\right) \right)
u^{2}\varphi _{\lambda ,\nu _{0}}\right) \right] + \\ 
+\partial _{t}\left( \left( \lambda ^{2}/\left( 2\lambda +1\right) \right)
u^{2}\varphi _{\lambda ,\nu _{0}}+\left( \nabla u\right) ^{2}\varphi
_{\lambda ,\nu _{0}}/\left( 2\lambda +1\right) \right) .%
\end{array}%
\right.  \label{4.30}
\end{equation}%
\emph{And the function }$\func{div}U$\emph{\ is:}%
\begin{equation}
\left. 
\begin{array}{c}
\func{div}U=\left[ \left( 2\lambda /\left( 2\lambda +1\right) \right) \left(
-2u_{t}\left( u_{x_{1}}+\lambda \nu \psi ^{\nu _{0}-1}u\right) \varphi
_{\lambda ,\nu _{0}}\psi ^{-\nu _{0}+1}\right) \right] _{x_{1}} \\ 
+\left[ \left( 2\lambda /\left( 2\lambda +1\right) \right) \left( -2\lambda
\nu _{0}\left( u_{x_{1}}+\lambda \nu _{0}\psi ^{\nu _{0}-1}u\right)
^{2}\varphi _{\lambda ,\nu _{0}}+2\lambda \nu
_{0}\dsum\limits_{i=2}^{n}u_{x_{i}}^{2}\varphi _{\lambda ,\nu _{0}}\right) %
\right] _{x_{1}}+ \\ 
+\left[ \left( 2\lambda /\left( 2\lambda +1\right) \right) \left( -2\lambda
^{3}\nu _{0}^{3}\psi ^{2\nu _{0}-2}\left( 1-2\psi ^{-\nu _{0}}\left( \nu
_{0}-1\right) /\left( \lambda \nu _{0}\right) \right) u^{2}\varphi _{\lambda
,\nu _{0}}\right) \right] _{x_{1}} \\ 
+\dsum\limits_{i=2}^{n}\left[ \left( 2\lambda /\left( 2\lambda +1\right)
\right) \left( -4\lambda \nu _{0}\left( u_{x_{1}}+\lambda \nu _{0}\psi ^{\nu
_{0}-1}u\right) u_{x_{i}}\varphi _{\lambda ,\nu _{0}}-2u_{t}u_{x_{i}}\varphi
_{\lambda ,\nu _{0}}\psi ^{-\nu _{0}+1}\right) \right] _{x_{i}}+ \\ 
+\left[ \left( 2\lambda /\left( 2\lambda +1\right) \right) \left( -\lambda
u_{x_{1}}u\varphi _{\lambda ,\nu _{0}}+\lambda ^{2}\nu _{0}\psi ^{\nu
_{0}-1}u^{2}\varphi _{\lambda ,\nu _{0}}\right) \right] _{x_{1}}+ \\ 
+\dsum\limits_{i=2}^{n}\left[ \left( 2\lambda /\left( 2\lambda +1\right)
\right) \left( -\lambda u_{x_{i}}u\varphi _{\lambda ,\nu _{0}}\right) \right]
_{x_{i}}+ \\ 
+\dsum\limits_{i=1}^{n}\left[ \left( 1/\left( 2\lambda +1\right) \right)
\left( -2u_{t}u_{x_{i}}\varphi _{\lambda ,\nu _{0}}\right) \right] _{x_{i}}+
\\ 
+\dsum\limits_{i=2}^{n}\left[ \left( 1/\left( 2\lambda +1\right) \right)
\left( -2u_{x_{1}x_{i}}u_{x_{i}}\varphi _{\lambda ,\nu _{0}}\right) \right]
_{x_{1}}+ \\ 
+\dsum\limits_{i=2}^{n}\left[ \left( 1/\left( 2\lambda +1\right) \right)
\left( 2u_{x_{1}x_{1}}u_{x_{i}}\varphi _{\lambda ,\nu _{0}}\right) \right]
_{x_{i}}+ \\ 
+\dsum\limits_{i,j=2}^{n}\left[ \left( 1/\left( 2\lambda +1\right) \right)
\left( u_{x_{j}x_{j}}u_{x_{i}}\varphi _{\lambda ,\nu
_{0}}-u_{x_{i}x_{j}}u_{x_{j}}\varphi _{\lambda ,\nu _{0}}\right) \right]
_{x_{i}}.%
\end{array}%
\right.  \label{4.31}
\end{equation}

\emph{In particular, (\ref{4.30}) leads to the following implications:}%
\begin{equation}
u\left( x,0\right) =u\left( x,T\right) \rightarrow V\left( x,0\right)
=V\left( x,T\right) \rightarrow \dint\limits_{Q_{T}}\partial _{t}Vdxdt=0.
\label{4.5}
\end{equation}

Below $C=C\left( Q_{T}\right) >0$ denotes different constants depending only
on the domain $\Omega .$ parameters.

\textbf{Remarks 4.1}:

\begin{enumerate}
\item \emph{Formula (\ref{4.30}) for the function }$V$ \emph{implies the key
property, which we need: that parasitic integrals over }$\left\{ t=0\right\} 
$\emph{\ and }$\left\{ t=T\right\} ,$\emph{\ which occur when integrating (%
\ref{4.3}) over }$Q_{T},$\emph{\ cancel each other, if }$u\left( x,T\right)
=u\left( x,0\right) ,$ \emph{and this is reflected in (\ref{4.5}) and in the
last line of (\ref{4.6}) of Theorem 4.2. The necessity of (\ref{4.5}) for
our goal of obtaining stability estimates for our CIP is explained in
subsection 3.2.}

\item \emph{Item 1 explains the reason of our need of a painstaking
derivation of the precise formula (\ref{4.30}) for }$\partial _{t}V$\emph{\
in the proof of Theorem 4.1. The reason of the derivation of precise formula
(\ref{4.31}) for }$\func{div}U$ \emph{is the necessity of the incorporation
of estimates of boundary terms, especially those with }$u_{t}$\emph{\ and }$%
u_{x_{i}x_{j}},$\emph{\ in the integral Carleman estimate (\ref{4.6}) of
Theorem 4.2.}
\end{enumerate}

\textbf{Theorem 4.2} (integral Carleman estimate for the operator $\partial
_{t}-\Delta ).$ \emph{Let }$\nu _{0}$\emph{\ and }$\lambda _{0}$\emph{\ be} 
\emph{parameters chosen }in \emph{Theorem 4.1. Then the following integral
Carleman estimate holds: }%
\begin{equation}
\left. 
\begin{array}{c}
C\exp \left[ 3\cdot 2^{\nu _{0}}\lambda \right] \left( \left\Vert
u\right\Vert _{H^{2,1}\left( \Gamma _{T}^{-}\right) }^{2}+\left\Vert
u_{x_{1}}\right\Vert _{H^{1,0}\left( \Gamma _{T}^{-}\right) }^{2}\right) +
\\ 
+C\exp \left[ 3\lambda \left( 2A_{1}+2\right) ^{\nu _{0}}\right] \left(
\left\Vert u\right\Vert _{H^{2,1}\left( S_{T}\diagdown \Gamma
_{T}^{-}\right) }^{2}+\left\Vert \partial _{n}u\right\Vert _{H^{1,0}\left(
S_{T}\diagdown \Gamma _{T}^{-}\right) }^{2}\right) + \\ 
+\dint\limits_{Q_{T}}\left( u_{t}-\Delta u\right) ^{2}\varphi _{\lambda ,\nu
_{0}}dxdt\geq \\ 
\geq \left( C/\lambda \right) \dint\limits_{Q_{T}}\left(
u_{t}^{2}+\dsum\limits_{i,j=2}^{n}u_{ij}^{2}\right) \varphi _{\lambda ,\nu
_{0}}dxdt+C\dint\limits_{Q_{T}}\left( \lambda \left( \nabla u\right)
^{2}+\lambda ^{3}u^{2}\right) \varphi _{\lambda ,\nu _{0}}dxdt, \\ 
\forall u\in H^{4,2}\left( Q_{T}\right) \cap \left\{ u:u\left( x,0\right)
=u\left( x,T\right) \right\} ,\forall \lambda \geq \lambda _{0}\emph{.}%
\end{array}%
\right.  \label{4.6}
\end{equation}

Since we have two parabolic operators in (\ref{2.2}), whose principal parts
are $\partial _{t}-\Delta $ and $\partial _{t}+\Delta ,$ then we need to
formulate an analog of Carleman estimate (\ref{4.6}) for the operator $%
\partial _{t}+\Delta $ as well. This is done in Theorem 4.3. We omit the
proof of this theorem, since it is quite similar with the proofs of Theorem
4.1 and 4.2. As to the norms involved in (\ref{4.6}), we refer to (\ref%
{2.101})-(\ref{2.105}).

\textbf{Theorem 4.3 }(integral Carleman estimates for the operator $\partial
_{t}+\Delta $). \emph{Let }$\nu _{0}$\emph{\ and }$\lambda _{0}$\emph{\ be
two parameters chosen in Theorem 4.1. \ Then the direct analog of the
Carleman estimate (\ref{4.6}) holds true when }$\left( u_{t}-\Delta u\right)
^{2}$\emph{\ is replaced with }$\left( u_{t}+\Delta u\right) ^{2}.$

\textbf{Remark 4.2.} \emph{We prove Theorems 4.1 and 4.2 in Appendix.
However, when carrying out other proofs below, we assume that Theorems 4.1
and 4.2 hold true.}

\subsection{H\"{o}lder and Lipschitz Stability estimates}

\label{sec:4.2}

In the theory of Ill-Posed Problems, one often assumes that solution of such
a problem belongs to an a priory chosen boundary set. Hence, let $M>0$ be
the number of subsection 2.4. Recalling Remark 1.1 and (\ref{2.100}), we
introduce the following set of pairs of functions $\left( u,m\right) :$%
\begin{equation}
Y_{1}\left( M\right) =\left\{ \left( u,m\right) \in C^{6,3}\left( \overline{Q%
}_{T}\right) :\left\Vert u\right\Vert _{C^{6,3}\left( \overline{Q}%
_{T}\right) },\left\Vert m\right\Vert _{C^{6,3}\left( \overline{Q}%
_{T}\right) }<M\right\} .  \label{4.7}
\end{equation}%
Obviously,%
\begin{equation}
\left\Vert u\right\Vert _{H^{6,3}\left( \overline{Q}_{T}\right) },\left\Vert
m\right\Vert _{H^{6,3}\left( \overline{Q}_{T}\right) }\leq CM,\text{ }%
\forall \left( u,m\right) \in Y_{1}\left( M\right) .  \label{4.8}
\end{equation}%
We also assume that functions $a\left( x\right) $ and $s\left( x,t\right) $
in MFGS (\ref{2.2}) satisfy the following conditions:%
\begin{equation}
a\in C^{3}\left( \overline{\Omega }\right) ,s\in C^{2,1}\left( \overline{Q}%
_{T}\right) ,\text{ }\left\Vert a\right\Vert _{C^{3}\left( \overline{\Omega }%
\right) }<M,\text{ }\left\Vert s\right\Vert _{C^{2,1}\left( \overline{Q}%
_{T}\right) }<M.  \label{4.9}
\end{equation}%
In addition, let the unknown coefficient%
\begin{equation}
b\left( x\right) \in Y_{2}\left( M\right) =\left\{ b:b\in C^{4}\left( 
\overline{\Omega }\right) ,\text{ }\left\Vert b\right\Vert _{C^{4}\left( 
\overline{\Omega }\right) }<M\right\} .  \label{4.10}
\end{equation}

\textbf{Theorem 4.4} (H\"{o}lder stability for incomplete data, the case (%
\ref{2.5})).\emph{\ Assume that there exists two vector functions }$\left(
u_{i},m_{i},b_{i}\right) \in Y_{1}\left( M\right) \times Y_{2}\left(
M\right) ,$\emph{\ }$i=1,2$\emph{\ satisfying the following analogs of
conditions (\ref{2.4}): }%
\begin{equation}
\left. 
\begin{array}{c}
u_{i}\left( x,0\right) =p_{i}\left( x\right) ,\text{ }m_{i}\left( x,0\right)
=q_{i}\left( x\right) ,\text{ }x\in \Omega , \\ 
u_{i}\left( x,T\right) =F_{i}\left( x\right) ,\text{ }m_{i}\left( x,T\right)
=G_{i}\left( x\right) ,\text{ }x\in \Omega , \\ 
u_{i}\mid _{S_{T}}=f_{0,i}\left( x,t\right) ,\text{ }\partial _{n}u_{i}\mid
_{S_{T}}=f_{1,i}\left( x,t\right) , \\ 
m_{i}\mid _{S_{T}}=g_{0,i}\left( x,t\right) ,\text{ }\partial _{n}m_{i}\mid
_{S_{T}}=g_{1,i}\left( x,t\right) , \\ 
i=1,2.%
\end{array}%
\right.  \label{4.11}
\end{equation}%
\emph{Assume that the lateral Cauchy data are incomplete as in the case (\ref%
{2.5}), i.e. functions }$f_{0i},f_{1i}$\emph{\ in (\ref{4.11}) are known for 
}$\left( x,t\right) \in S_{T}\diagdown \Gamma _{T}^{-}$\emph{\ and are
unknown for }$\left( x,t\right) \in \Gamma _{T}^{-}.$\emph{\ Let }$\delta
\in \left( 0,1\right) $\emph{\ be a number characterizing the level of the
error in the data} (\ref{4.11}). \emph{More precisely, let}%
\begin{equation}
\left. 
\begin{array}{c}
\left\Vert p_{1}-p_{2}\right\Vert _{H^{4}\left( \Omega \right) }<\delta
,\left\Vert q_{1}-q_{2}\right\Vert _{H^{3}\left( \Omega \right) }<\delta ,
\\ 
\left\Vert F_{1}-F_{2}\right\Vert _{H^{4}\left( \Omega \right) }<\delta
,\left\Vert G_{1}-G_{2}\right\Vert _{H^{3}\left( \Omega \right) }<\delta ,
\\ 
\left\Vert \partial _{t}f_{0,1}-\partial _{t}f_{0,2}\right\Vert
_{H^{2,1}\left( S_{T}\diagdown \Gamma _{T}^{-}\right) }<\delta ,\left\Vert
\partial _{t}f_{1,1}-\partial _{t}f_{1,2}\right\Vert _{H^{1,0}\left(
S_{T}\diagdown \Gamma _{T}^{-}\right) }<\delta , \\ 
\left\Vert \partial _{t}g_{0,1}-\partial _{t}g_{0,2}\right\Vert
_{H^{2,1}\left( S_{T}\diagdown \Gamma _{T}^{-}\right) }<\delta ,\left\Vert
\partial _{t}g_{1,1}-\partial _{t}g_{1,2}\right\Vert _{H^{1,0}\left(
S_{T}\diagdown \Gamma _{T}^{-}\right) }<\delta .%
\end{array}%
\right.  \label{4.12}
\end{equation}%
\emph{Assume that condition (\ref{2.3}) holds. In addition, assume that
there exists a number }$c>0$\emph{\ such that} 
\begin{equation}
\min_{\overline{Q}_{T}}\left\vert \dint\limits_{\Omega ^{\prime
}}K_{1}\left( \overline{x},\overline{y}\right) m_{2}\left( x_{1},\overline{y}%
,t\right) d\overline{y}+\dint\limits_{x_{1}}^{A_{1}}\left(
\dint\limits_{\Omega ^{\prime }}K_{2}\left( x,y_{1},\overline{y}\right)
m_{2}\left( y_{1},\overline{y},t\right) d\overline{y}\right)
dy_{1}\right\vert \geq c.  \label{4.13}
\end{equation}%
\emph{Let }$\gamma \in \left( 0,2A_{1}\right) $\emph{\ be the number in (\ref%
{2.1}). Then there exist a sufficiently small number }$\delta _{0}=\delta
_{0}\left( M,c,\gamma ,\Omega ,T\right) \in \left( 0,1\right) $\emph{\ and a
number }$B=B\left( M,c,\gamma ,\Omega ,T\right) >0,$\emph{\ both numbers
depending only on listed parameters, such that for all }$\delta \in \left(
0,\delta _{0}\right) $\emph{\ the following H\"{o}lder stability estimates
are valid with a certain number }$\alpha \in \left( 0,1\right) :$%
\begin{equation}
\left\Vert \partial _{t}^{j}u_{1}-\partial _{t}^{j}u_{2}\right\Vert
_{H^{2,1}\left( Q_{\gamma T}\right) },\left\Vert \partial
_{t}^{j}m_{1}-\partial _{t}^{j}m_{2}\right\Vert _{H^{2,1}\left( Q_{\gamma
T}\right) }\leq B\delta ^{\alpha },\text{ }j=0,1,  \label{4.14}
\end{equation}%
\begin{equation}
\left\Vert b_{1}-b_{2}\right\Vert _{L_{2}\left( \Omega _{\gamma }\right)
}\leq B\delta ^{\alpha }.  \label{4.15}
\end{equation}%
\emph{In particular, our CIP with the incomplete data as in (\ref{2.1}) has
at most one solution.}

\textbf{Remark 4.3.} \emph{Below }$B=B\left( M,c,\gamma ,\Omega ,T\right) >0$%
\emph{\ and }$C_{1}=C_{1}\left( M,c,\Omega ,T\right) >0$\emph{\ denote
different numbers depending only on listed parameters.}

\textbf{Theorem 4.5} (Lipschitz stability for complete data, the case (\ref%
{2.6})). \emph{Assume that there exists two vector functions }$\left(
u_{i},m_{i},b_{i}\right) \in Y_{1}\left( M\right) \times Y_{2}\left(
M\right) ,$\emph{\ }$i=1,2$\emph{\ satisfying conditions (\ref{4.11}).
Assume that the lateral Cauchy data are complete as in (\ref{2.6}), i.e.
functions }$f_{0i},f_{1i},g_{0i},g_{1i}$\emph{\ in (\ref{4.11}) are known
for all }$\left( x,t\right) \in S_{T}.$\emph{\ In addition, let conditions (%
\ref{2.3}) and (\ref{4.13}) hold. Then the following Lipschitz stability
estimates are valid:}%
\begin{equation*}
\left. 
\begin{array}{c}
\left\Vert \partial _{t}^{j}u_{1}-\partial _{t}^{j}u_{2}\right\Vert
_{H^{2,1}\left( Q_{T}\right) },\left\Vert \partial _{t}^{j}m_{1}-\partial
_{t}^{j}m_{2}\right\Vert _{H^{2,1}\left( Q_{T}\right) }\leq \\ 
\leq C_{1}\left( \left\Vert p_{1}-p_{2}\right\Vert _{H^{4}\left( \Omega
\right) }+\left\Vert q_{1}-q_{2}\right\Vert _{H^{4}\left( \Omega \right)
}\right) + \\ 
+C_{1}\left( \left\Vert F_{1}-F_{2}\right\Vert _{H^{4}\left( \Omega \right)
}+\left\Vert G_{1}-G_{2}\right\Vert _{H^{4}\left( \Omega \right) }\right) +
\\ 
+C_{1}\left( \left\Vert \partial _{t}f_{0,1}-\partial _{t}f_{0,2}\right\Vert
_{H^{2,1}\left( S_{T}{}^{-}\right) }+\left\Vert \partial
_{t}f_{1,1}-\partial _{t}f_{1,2}\right\Vert _{H^{1,0}\left( S_{T}\right)
}\right) + \\ 
+C_{1}\left( \left\Vert \partial _{t}g_{0,1}-\partial _{t}g_{0,2}\right\Vert
_{H^{2,1}\left( S_{T}{}^{-}\right) }+\left\Vert \partial
_{t}g_{1,1}-\partial _{t}g_{1,2}\right\Vert _{H^{1,0}\left( S_{T}\right)
}\right) , \\ 
j=0,1.%
\end{array}%
\right.
\end{equation*}%
\begin{equation*}
\left. 
\begin{array}{c}
\left\Vert b_{1}-b_{2}\right\Vert _{L_{2}\left( \Omega \right) }\leq \\ 
\leq \left\Vert p_{1}-p_{2}\right\Vert _{H^{4}\left( \Omega \right)
}+\left\Vert q_{1}-q_{2}\right\Vert _{H^{4}\left( \Omega \right) }+ \\ 
+\left\Vert F_{1}-F_{2}\right\Vert _{H^{4}\left( \Omega \right) }+\left\Vert
G_{1}-G_{2}\right\Vert _{H^{4}\left( \Omega \right) }+ \\ 
+\left\Vert \partial _{t}f_{0,1}-\partial _{t}f_{0,2}\right\Vert
_{H^{2,1}\left( S_{T}{}^{-}\right) }+\left\Vert \partial
_{t}f_{1,1}-\partial _{t}f_{1,2}\right\Vert _{H^{1,0}\left( S_{T}\right) }.%
\end{array}%
\right.
\end{equation*}

\section{Proofs of Theorems 4.4 and 4.5}

\label{sec:5}

\subsection{Proof of Theorem 4.4}

\QTP{Body Math}
\label{sec:5.1}

First, we proceed with the transformation procedure as outlined in
subsection 3.2. Next, we apply to two resulting integral differential
equations Carleman estimates of Theorems 4.2 and 4.3.

\subsubsection{The transformation procedure}

\label{sec:5.1.1}

Consider the differences:%
\begin{equation}
\left. 
\begin{array}{c}
\widetilde{u}=u_{1}-u_{2},\text{ }\widetilde{m}=m_{1}-m_{2},\text{ }%
\widetilde{b}=b_{1}-b_{2}, \\ 
\widetilde{p}=p_{1}-p_{2},\text{ }\widetilde{q}=q_{1}-q_{2}, \\ 
\widetilde{F}=F_{1}-F_{2},\text{ }\widetilde{G}=G_{1}-G_{2}, \\ 
\widetilde{f}_{0}=f_{0,1}-f_{0,2},\text{ }\widetilde{f}_{1}=f_{1,1}-f_{1,2},
\\ 
\widetilde{g}_{0}=g_{0,1}-g_{0,2},\text{ }\widetilde{g}_{1}=g_{1,1}-g_{1,2}.%
\end{array}%
\right.  \label{5.1}
\end{equation}%
By (\ref{2.5}), (\ref{4.11}) and (\ref{5.1})%
\begin{equation}
\left. 
\begin{array}{c}
\widetilde{u}\left( x,0\right) =\widetilde{p}\left( x\right) ,\text{ }%
\widetilde{m}\left( x,0\right) =\widetilde{q}\left( x\right) ,\text{ }x\in
\Omega , \\ 
\widetilde{u}\left( x,T\right) =\widetilde{F}\left( x\right) ,\text{ }%
\widetilde{m}\left( x,T\right) =\widetilde{G}\left( x\right) ,\text{ }x\in
\Omega , \\ 
\widetilde{u}\mid _{S_{T}\diagdown \Gamma _{T}^{-}}=\widetilde{f}_{0}\left(
x,t\right) ,\text{ }\partial _{n}\widetilde{u}\mid _{S_{T}\diagdown \Gamma
_{T}^{-}}=\widetilde{f}_{1}\left( x,t\right) , \\ 
\widetilde{m}\mid _{S_{T}\diagdown \Gamma _{T}^{-}}=\widetilde{g}_{0}\left(
x,t\right) ,\text{ }\partial _{n}\widetilde{m}\mid _{S_{T}\diagdown \Gamma
_{T}^{-}}=\widetilde{g}_{1}\left( x,t\right) .%
\end{array}%
\right.  \label{5.2}
\end{equation}

Let $y_{1},z_{1}$ and $y_{2},z_{2}$ be two pairs of numbers. Denote $%
\widetilde{y}=y_{1}-y_{2},\widetilde{z}=z_{1}-z_{2}.$ Then%
\begin{equation}
y_{1}z_{1}-y_{2}z_{2}=\widetilde{y}z_{1}+\widetilde{z}y_{2}.  \label{5.3}
\end{equation}%
Subtracting equations (\ref{2.2}) for $\left( u_{2},m_{2},b_{2}\right) $
from the same equations for $\left( u_{1},m_{1},b_{1}\right) $ and using (%
\ref{2.30}), the first line of (\ref{5.1}) as well as (\ref{5.3}), we obtain
two equations (\ref{5.4}) and (\ref{5.5}), 
\begin{equation}
\left. 
\begin{array}{c}
\widetilde{u}_{t}(x,t)+\Delta \widetilde{u}(x,t){-a(x)\nabla }\left(
u_{1}+u_{2}\right) {(\nabla \widetilde{u}(x,t)/2}+ \\ 
+b_{1}\left( x\right) \dint\limits_{\Omega ^{\prime }}K_{1}\left( \overline{x%
},\overline{y}\right) \widetilde{m}\left( x_{1},\overline{y},t\right) d%
\overline{y}+ \\ 
+b_{1}\left( x\right) \dint\limits_{x_{1}}^{A_{1}}\left(
\dint\limits_{\Omega ^{\prime }}K_{2}\left( x,y_{1},\overline{y}\right) 
\widetilde{m}\left( y_{1},\overline{y},t\right) d\overline{y}\right) dy_{1}=
\\ 
=-\widetilde{b}\left( x\right) \dint\limits_{\Omega ^{\prime }}K_{1}\left( 
\overline{x},\overline{y}\right) m_{2}\left( x_{1},\overline{y},t\right) d%
\overline{y}- \\ 
-\widetilde{b}\left( x\right) \dint\limits_{x_{1}}^{A_{1}}\left(
\dint\limits_{\Omega ^{\prime }}K_{2}\left( x,y_{1},\overline{y}\right)
m_{2}\left( y_{1},\overline{y},t\right) d\overline{y}\right) dy_{1}, \\ 
\left( x,t\right) \in Q_{T},%
\end{array}%
\right.  \label{5.4}
\end{equation}%
\begin{equation}
\begin{array}{c}
\widetilde{m}_{t}(x,t)-\Delta \widetilde{m}(x,t){-\func{div}(a(x)\widetilde{m%
}(x,t)\nabla u}_{2}{(x,t))-} \\ 
-{\func{div}(a(x)m}_{1}{(x,t)\nabla \widetilde{u}(x,t))}=0,\text{ } \\ 
\left( x,t\right) \in Q_{T}.%
\end{array}
\label{5.5}
\end{equation}%
Divide both sides of equation (\ref{5.4}) by the function $R\left(
x,t\right) ,$%
\begin{equation}
R\left( x,t\right) =-\dint\limits_{\Omega ^{\prime }}K_{1}\left( \overline{x}%
,\overline{y}\right) m_{2}\left( x_{1},\overline{y},t\right) d\overline{y}%
-\dint\limits_{x_{1}}^{A_{1}}\left( \dint\limits_{\Omega ^{\prime
}}K_{2}\left( x,y_{1},\overline{y}\right) m_{2}\left( y_{1},\overline{y}%
,t\right) d\overline{y}\right) dy_{1}.  \label{5.6}
\end{equation}%
By (\ref{4.13}) and (\ref{5.6}) 
\begin{equation}
\frac{1}{\left\vert R\left( x,t\right) \right\vert }\geq \frac{1}{c}.
\label{5.7}
\end{equation}%
Denote%
\begin{equation}
\overline{u}\left( x,t\right) =\frac{\widetilde{u}\left( x,t\right) }{%
R\left( x,t\right) }.  \label{5.8}
\end{equation}%
Then equation (\ref{5.4}) becomes:%
\begin{equation}
\left. 
\begin{array}{c}
\overline{u}_{t}+\Delta \overline{u}+P\nabla \overline{u}+Q\overline{u}+ \\ 
+b_{1}\left( x\right) R^{-1}\left( x,t\right) \dint\limits_{\Omega ^{\prime
}}K_{1}\left( \overline{x},\overline{y}\right) \widetilde{m}\left( x_{1},%
\overline{y},t\right) d\overline{y}+ \\ 
+b_{1}\left( x\right) R^{-1}\left( x,t\right)
\dint\limits_{x_{1}}^{A_{1}}\left( \dint\limits_{\Omega ^{\prime
}}K_{2}\left( x,y_{1},\overline{y}\right) \widetilde{m}\left( y_{1},%
\overline{y},t\right) d\overline{y}\right) dy_{1}= \\ 
=-\widetilde{b}\left( x\right) ,%
\end{array}%
\right.  \label{5.9}
\end{equation}%
where $\left( P,Q\right) $ is an $\left( n+1\right) -$dimensional vector
function with its $C^{2,1}\left( \overline{Q}_{T}\right) -$components.
Although it is easy to present the explicit formulas for its components, we
are not doing this for brevity.

It follows from (\ref{5.2}), (\ref{5.8}) and (\ref{5.9}) that the function $%
\overline{u}_{t}\left( x,t\right) $ attains the following values at $t=0,T:$%
\begin{equation}
\left. 
\begin{array}{c}
\overline{u}_{t}\left( x,0\right) = \\ 
=-\Delta \left( R^{-1}\left( x,0\right) \widetilde{p}\left( x\right) \right)
-P\left( x,0\right) \nabla \left( R^{-1}\left( x,0\right) \widetilde{p}%
\left( x\right) \right) - \\ 
-Q\left( x,0\right) \nabla \left( R^{-1}\left( x,0\right) \widetilde{p}%
\left( x\right) \right) -b_{1}\left( x\right) R^{-1}\left( x,0\right)
\dint\limits_{\Omega ^{\prime }}K_{1}\left( \overline{x},\overline{y}\right) 
\widetilde{q}\left( x_{1},\overline{y}\right) d\overline{y}- \\ 
-b_{1}\left( x\right) R^{-1}\left( x,0\right)
\dint\limits_{x_{1}}^{A_{1}}\left( \dint\limits_{\Omega ^{\prime
}}K_{2}\left( x,y_{1},\overline{y}\right) \widetilde{q}\left( y_{1},%
\overline{y}\right) d\overline{y}\right) dy_{1}-\widetilde{b}\left( x\right)
, \\ 
=W_{0}\left( x\right) -\widetilde{b}\left( x\right) , \\ 
\begin{array}{c}
\overline{u}_{t}\left( x,T\right) =-\Delta \left( R^{-1}\left( x,T\right) 
\widetilde{F}\left( x\right) \right) -P\left( x,T\right) \nabla \left(
R^{-1}\left( x,0\right) \widetilde{F}\left( x\right) \right) - \\ 
-Q\left( x,T\right) \nabla \left( R^{-1}\left( x,T\right) \widetilde{F}%
\left( x\right) \right) -b_{1}\left( x\right) R^{-1}\left( x,T\right)
\dint\limits_{\Omega ^{\prime }}K_{1}\left( \overline{x},\overline{y}\right) 
\widetilde{G}\left( x_{1},\overline{y}\right) d\overline{y}- \\ 
-b_{1}\left( x\right) R^{-1}\left( x,T\right)
\dint\limits_{x_{1}}^{A_{1}}\left( \dint\limits_{\Omega ^{\prime
}}K_{2}\left( x,y_{1},\overline{y}\right) \widetilde{G}\left( y_{1},%
\overline{y}\right) d\overline{y}\right) dy_{1}-\widetilde{b}\left( x\right)
= \\ 
=W_{T}\left( x\right) -\widetilde{b}\left( x\right) ,%
\end{array}%
\end{array}%
\right.  \label{5.10}
\end{equation}%
Next, (\ref{4.11}), (\ref{5.2}), (\ref{5.5}) and (\ref{5.8}) imply similar
formulas for ${\widetilde{m}}_{t}{(x,0)}$ and ${\widetilde{m}}_{t}{(x,T),}$%
\begin{equation}
\begin{array}{c}
\widetilde{m}_{t}(x,0)= \\ 
=\Delta \widetilde{q}(x)+{\func{div}(a(x)\widetilde{q}(x)\nabla p}_{2}{(x))-}
\\ 
+{\func{div}(a(x)q}_{1}\left( x\right) {\nabla }\left( \widetilde{p}\left(
x\right) /R\left( x,0\right) \right) = \\ 
=Z_{0}\left( x\right) ,\text{ } \\ 
\\ 
\begin{array}{c}
\widetilde{m}_{t}(x,T)= \\ 
=\Delta \widetilde{G}(x)+{\func{div}(a(x)\widetilde{G}(x)\nabla F}_{2}{(x))-}
\\ 
+{\func{div}(a(x)G}_{1}\left( x\right) {\nabla }\left( \widetilde{F}\left(
x\right) /R\left( x,T\right) \right) = \\ 
=Z_{T}\left( x\right) .%
\end{array}%
\end{array}
\label{5.11}
\end{equation}

Differentiate both equations (\ref{5.5}) and (\ref{5.9}) with respect to $t.$
Then we obtain equations for the $t-$derivatives $\overline{v}\left(
x,t\right) $ and $\overline{w}\left( x,t\right) .$ An important property of
these equations is that the function $\widetilde{b}\left( x\right) $ is not
present in them, since this function is independent on $t$. Functions $%
\overline{w}\left( x,t\right) $ and $\overline{w}\left( x,t\right) $ are: 
\begin{equation}
\text{ }\overline{v}\left( x,t\right) =\overline{u}_{t}\left( x,t\right) ,%
\overline{w}\left( x,t\right) =\widetilde{m}_{t}\left( x,t\right)
\label{5.12}
\end{equation}

By the first line of (\ref{5.2}) as well as by (\ref{5.8}) and (\ref{5.12})%
\begin{equation}
\text{ }\overline{u}\left( x,t\right) =\dint\limits_{0}^{t}\overline{v}%
\left( x,\tau \right) d\tau +\frac{\widetilde{p}\left( x\right) }{R\left(
x,0\right) },\text{ }\widetilde{m}\left( x,t\right) =\dint\limits_{0}^{t}%
\overline{w}\left( x,\tau \right) d\tau +\widetilde{q}\left( x\right) .
\label{5.13}
\end{equation}%
Substitute (\ref{5.13}) in equations for functions $\overline{v}\left(
x,t\right) $ and $\overline{w}\left( x,t\right) .$ Then introduce new
functions%
\begin{equation}
v\left( x,t\right) =\overline{v}\left( x,t\right) -\left( W_{0}\left(
x\right) \frac{t}{T}+W_{T}\left( x\right) \left( 1-\frac{t}{T}\right)
\right) ,  \label{5.14}
\end{equation}%
\begin{equation}
w\left( x,t\right) =\overline{w}\left( x,t\right) -\left( Z_{0}\left(
x\right) \frac{t}{T}+Z_{T}\left( x\right) \left( 1-\frac{t}{T}\right)
\right) ,  \label{5.15}
\end{equation}%
where functions $W_{0}\left( x\right) ,W_{T}\left( x\right) ,Z_{0}\left(
x\right) ,Z_{T}\left( x\right) $ are given in (\ref{5.10}) and (\ref{5.11}).
It follows from (\ref{5.10}), (\ref{5.11}), (\ref{5.14}) and (\ref{5.15})
that 
\begin{equation}
\left. 
\begin{array}{c}
v\left( x,0\right) =v\left( x,T\right) =-\widetilde{b}\left( x\right) , \\ 
w\left( x,0\right) =w\left( x,T\right) =0.%
\end{array}%
\right.  \label{5.16}
\end{equation}%
This finishes our transformation procedure outlined in subsection 3.2.
Indeed, comparing the last line of (\ref{4.6}) with (\ref{5.16}), we see
that Carleman estimates of Theorems 4.2 and 4.3 can be applied to functions $%
w$ and $v$ respectively.

\paragraph{Applications of Theorems 4.2 and 4.3}

\label{sec:5.1.2}

It is well known that Carleman estimates can work not only with equations
but with inequalities as well. Hence, to simplify the presentation, we turn
the above mentioned equations for functions $v$ and $w$ in two integral
differential inequalities. The first inequality is:%
\begin{equation}
\left. 
\begin{array}{c}
\left\vert v_{t}+\Delta v\right\vert \leq C_{1}\left( \left\vert \nabla
v\right\vert +\left\vert v\right\vert \right) + \\ 
+C_{1}\dint\limits_{\Omega ^{\prime }}\left\vert w\left( x_{1},\overline{y}%
,t\right) \right\vert d\overline{y}+C_{1}\dint\limits_{\Omega ^{\prime
}}\left( \dint\limits_{0}^{t}\left\vert w\left( x_{1},\overline{y},\tau
\right) \right\vert d\tau \right) d\overline{y}+ \\ 
+C_{1}\dint\limits_{x_{1}}^{A_{1}}\left( \dint\limits_{\Omega ^{\prime
}}\left\vert w\right\vert \left( y_{1},\overline{y},t\right) d\overline{y}%
\right) dy_{1}+C_{1}\dint\limits_{x_{1}}^{A_{1}}\left[ \dint\limits_{\Omega
^{\prime }}\left( \dint\limits_{0}^{t}\left\vert w\right\vert \left( y_{1},%
\overline{y},\tau \right) d\tau \right) d\overline{y}\right] dy_{1}+ \\ 
+X_{1}\left( x,t\right) ,\text{ }\left( x,t\right) \in Q_{T}.%
\end{array}%
\right.  \label{5.17}
\end{equation}

\bigskip The second inequality is:%
\begin{equation}
\left. 
\begin{array}{c}
\left\vert w_{t}-\Delta w\right\vert \leq C_{1}\left( \left\vert \nabla
w\right\vert +\left\vert w\right\vert +\dint\limits_{0}^{t}\left( \left\vert
\nabla w\right\vert +\left\vert w\right\vert \right) \left( x,\tau \right)
d\tau \right) + \\ 
+C_{1}\left( \left\vert \nabla v\right\vert +\dint\limits_{0}^{t}\left\vert
\nabla v\right\vert \left( x,\tau \right) d\tau \right) +C_{1}\left(
\left\vert \Delta v\right\vert +\dint\limits_{0}^{t}\left\vert \Delta
v\right\vert \left( x,\tau \right) d\tau \right) +X_{2}\left( x,t\right) ,
\\ 
\left( x,t\right) \in Q_{T}.%
\end{array}%
\right.  \label{5.18}
\end{equation}

Note that the presence of integrals with respect to $y$ in the third line of
(\ref{5.17}) is the indication of the third difficulty of working with MFGS (%
\ref{2.2}) as in item 3 of subsection 2.3. And the presence of the terms
with $\left\vert \Delta v\right\vert $ in (\ref{5.18}) is the indication of
the fourth difficulty listed in item 4 of that subsection.

It easily follows from the second and third lines of (\ref{2.3}), (\ref{4.7}%
)-(\ref{4.10}) and (\ref{5.6})-(\ref{5.15}) that functions $X_{1}$ and $%
X_{2} $ in (\ref{5.17}) and (\ref{5.18}) are such that 
\begin{equation}
\left. 
\begin{array}{c}
X_{1},X_{2}\in L_{2}\left( Q_{T}\right) , \\ 
\text{ }\left\Vert X_{1}\right\Vert _{L_{2}\left( Q_{T}\right)
}^{2}+\left\Vert X_{2}\right\Vert _{L_{2}\left( Q_{T}\right) }^{2}\leq \\ 
\leq C_{1}\left( \left\Vert \widetilde{p}\right\Vert _{H^{4}\left( \Omega
\right) }^{2}+\left\Vert \widetilde{q}\right\Vert _{H^{4}\left( \Omega
\right) }^{2}+\left\Vert \widetilde{F}\right\Vert _{H^{4}\left( \Omega
\right) }^{2}+\left\Vert \widetilde{G}\right\Vert _{H^{4}\left( \Omega
\right) }^{2}\right) .%
\end{array}%
\right.  \label{5.19}
\end{equation}

As to the lateral Cauchy data for functions $v$ and $w$, using boundary data
in (\ref{5.1}), (\ref{5.2}) and the above transformation procedure combined
with the considerations, which resulted in estimates (\ref{5.19}), we obtain%
\begin{equation}
\begin{array}{c}
v\mid _{S_{T}\diagdown \Gamma _{T}^{-}}=m_{0}\left( x,t\right) ,\text{ }%
\partial _{n}v\mid _{S_{T}\diagdown \Gamma _{T}^{-}}=m_{1}\left( x,t\right) ,
\\ 
w\mid _{S_{T}\diagdown \Gamma _{T}^{-}}=z_{0}\left( x,t\right) ,\text{ }%
\partial _{n}w\mid _{S_{T}\diagdown \Gamma _{T}^{-}}=z_{1}\left( x,t\right) ,%
\end{array}
\label{5.20}
\end{equation}%
where functions in the right hand sides of (\ref{5.20}) can be estimated as:%
\begin{equation}
\left. 
\begin{array}{c}
\left\Vert m_{0}\right\Vert _{H^{2,1}\left( S_{T}\diagdown \Gamma
_{T}^{-}\right) }^{2}+\left\Vert z_{0}\right\Vert _{H^{2,1}\left(
S_{T}\diagdown \Gamma _{T}^{-}\right) }^{2}+\left\Vert m_{1}\right\Vert
_{H^{1,0}\left( S_{T}\diagdown \Gamma _{T}^{-}\right) }^{2}+\left\Vert
z_{1}\right\Vert _{H^{1,0}\left( S_{T}\diagdown \Gamma _{T}^{-}\right)
}^{2}\leq \\ 
\leq C_{1}\left( \left\Vert \widetilde{f}_{0}\right\Vert _{H^{2,1}\left(
S_{T}\diagdown \Gamma _{T}^{-}\right) }^{2}+\left\Vert \widetilde{g}%
_{0}\right\Vert _{H^{2,1}\left( S_{T}\diagdown \Gamma _{T}^{-}\right)
}^{2}\right) + \\ 
+C_{1}\left( \left\Vert \widetilde{f}_{1}\right\Vert _{H^{1,0}\left(
S_{T}\diagdown \Gamma _{T}^{-}\right) }^{2}+\left\Vert \widetilde{g}%
_{1}\right\Vert _{H^{1,0}\left( S_{T}\diagdown \Gamma _{T}^{-}\right)
}^{2}\right) + \\ 
+C_{1}\left( \left\Vert \widetilde{p}\right\Vert _{H^{4}\left( \Omega
\right) }^{2}+\left\Vert \widetilde{q}\right\Vert _{H^{4}\left( \Omega
\right) }^{2}+\left\Vert \widetilde{F}\right\Vert _{H^{4}\left( \Omega
\right) }^{2}+\left\Vert \widetilde{G}\right\Vert _{H^{4}\left( \Omega
\right) }^{2}\right) .%
\end{array}%
\right.  \label{5.21}
\end{equation}

Square both sides of each of inequalities (\ref{5.17}) and (\ref{5.18}),
multiply by the CWF $\varphi _{\lambda ,\nu _{0}}\left( x\right) $ in (\ref%
{4.1}) and integrate over the domain $Q_{T}.$ Use Cauchy-Schwarz inequality,
the right inequality (\ref{4.2}) and (\ref{5.19}). Also, note that since the
function $\varphi _{\lambda ,\nu _{0}}\left( x\right) $ depends only on $%
x_{1},$ then%
\begin{equation}
\left. 
\begin{array}{c}
\dint\limits_{Q_{T}}\left( \dint\limits_{0}^{t}\left\vert f\right\vert
\left( x,\tau \right) d\tau \right) ^{2}\varphi _{\lambda ,\nu _{0}}\left(
x\right) dxdt\leq \\ 
\leq C_{1}\dint\limits_{Q_{T}}f^{2}\left( x,t\right) \varphi _{\lambda ,\nu
_{0}}\left( x\right) dxdt,\text{ }\forall f\in L_{2}\left( Q_{T}\right) , \\ 
\dint\limits_{Q_{T}}\left( \dint\limits_{\Omega ^{\prime }}\left\vert
f\left( x_{1},\overline{y},t\right) \right\vert d\overline{y}\right)
^{2}\varphi _{\lambda ,\nu _{0}}\left( x\right) dxdt\leq \\ 
\leq C_{1}\dint\limits_{Q_{T}}f^{2}\left( x,t\right) \varphi _{\lambda ,\nu
_{0}}\left( x\right) dxdt,\text{ }\forall f\in L_{2}\left( Q_{T}\right) .%
\end{array}%
\right.  \label{5.210}
\end{equation}%
Hence, we obtain two inequalities. The first inequality is:%
\begin{equation}
\left. 
\begin{array}{c}
\dint\limits_{Q_{T}}\left( v_{t}+\Delta v\right) ^{2}\varphi _{\lambda ,\nu
_{0}}\left( x\right) dxdt\leq C_{1}\dint\limits_{Q_{T}}\left( \left( \nabla
v\right) ^{2}+v^{2}\right) \varphi _{\lambda ,\nu _{0}}\left( x\right) dxdt+
\\ 
+C_{1}\dint\limits_{Q_{T}}w^{2}\varphi _{\lambda ,\nu _{0}}\left( x\right)
dxdt+C_{1}\dint\limits_{Q_{T}}\left(
\dint\limits_{x_{1}}^{A_{1}}\dint\limits_{\Omega ^{\prime }}w^{2}\left(
y_{1},\overline{y},t\right) d\overline{y}dy_{1}\right) \varphi _{\lambda
,\nu _{0}}\left( x\right) dxdt+ \\ 
+C_{1}\dint\limits_{Q_{T}}\left[ \dint\limits_{x_{1}}^{A_{1}}\dint\limits_{%
\Omega ^{\prime }}\left( \dint\limits_{0}^{t}w^{2}\left( y_{1},\overline{y}%
,\tau \right) d\tau \right) d\overline{y}dy_{1}\right] \varphi _{\lambda
,\nu _{0}}\left( x\right) dxdt \\ 
+C_{1}\left( \left\Vert \widetilde{p}\right\Vert _{H^{4}\left( Q_{T}\right)
}^{2}+\left\Vert \widetilde{q}\right\Vert _{H^{4}\left( Q_{T}\right)
}^{2}+\right) \exp \left[ 2\lambda \left( 2A_{1}+2\right) ^{\nu _{0}}\right]
+ \\ 
+C_{1}\left( \left\Vert \widetilde{F}\right\Vert _{H^{4}\left( Q_{T}\right)
}^{2}+\left\Vert \widetilde{G}\right\Vert _{H^{4}\left( Q_{T}\right)
}^{2}\right) \exp \left[ 2\lambda \left( 2A_{1}+2\right) ^{\nu _{0}}\right] .%
\end{array}%
\right.  \label{5.22}
\end{equation}%
We now estimate from the above the term in the third line of (\ref{5.22}).
When doing so, we recall that (\ref{4.1}) implies that the function $\varphi
_{\lambda ,\nu _{0}}\left( x\right) \equiv \varphi _{\lambda ,\nu
_{0}}\left( x_{1}\right) $ is increasing. Using (\ref{2.1}), we obtain%
\begin{equation*}
\left. 
\begin{array}{c}
\dint\limits_{Q_{T}}\left[ \dint\limits_{x_{1}}^{A_{1}}\dint\limits_{\Omega
^{\prime }}\left( \dint\limits_{0}^{t}w^{2}\left( y_{1},\overline{y},\tau
\right) d\tau \right) d\overline{y}dy_{1}\right] \varphi _{\lambda ,\nu
_{0}}\left( x\right) dxdt\leq \\ 
\leq T\dint\limits_{0|}^{T}dt\dint\limits_{\Omega ^{\prime }}d\overline{x}%
\dint\limits_{\Omega ^{\prime }}d\overline{y}\dint\limits_{-A_{1}}^{A_{1}}%
\varphi _{\lambda ,\nu _{0}}\left( x_{1}\right) \left(
\dint\limits_{x_{1}}^{A_{1}}w^{2}\left( y_{1},\overline{y},t\right)
dy_{1}\right) dx_{1}= \\ 
=T\dint\limits_{0|}^{T}dt\dint\limits_{\Omega ^{\prime }}d\overline{x}%
\dint\limits_{\Omega ^{\prime }}d\overline{y}\dint\limits_{-A_{1}}^{A_{1}}%
\left( \dint\limits_{-A_{1}}^{y_{1}}\varphi _{\lambda ,\nu _{0}}\left(
x_{1}\right) dx_{1}\right) w^{2}\left( y_{1},\overline{y},t\right) dy_{1}\leq
\\ 
\leq \dint\limits_{0|}^{T}dt\dint\limits_{\Omega ^{\prime }}d\overline{x}%
\dint\limits_{\Omega ^{\prime }}d\overline{y}\dint%
\limits_{-A_{1}}^{A_{1}}w^{2}\left( y_{1},\overline{y},t\right) \left(
y_{1}+A_{1}\right) \varphi _{\lambda ,\nu _{0}}\left( y_{1}\right) dy_{1}\leq
\\ 
\leq C_{1}\dint\limits_{Q_{T}}w^{2}\varphi _{\lambda ,\nu _{0}}\left(
x\right) dxdt.%
\end{array}%
\right.
\end{equation*}%
Hence, (\ref{5.22}) can be rewritten as:%
\begin{equation}
\left. 
\begin{array}{c}
\dint\limits_{Q_{T}}\left( v_{t}+\Delta v\right) ^{2}\varphi _{\lambda ,\nu
_{0}}\left( x\right) dxdt\leq C_{1}\dint\limits_{Q_{T}}\left( \left( \nabla
v\right) ^{2}+v^{2}\right) \varphi _{\lambda ,\nu _{0}}\left( x\right) dxdt+
\\ 
+C_{1}\dint\limits_{Q_{T}}w^{2}\varphi _{\lambda ,\nu _{0}}\left( x\right)
dxdt+ \\ 
+C_{1}\left( \left\Vert \widetilde{p}\right\Vert _{H^{4}\left( Q_{T}\right)
}^{2}+\left\Vert \widetilde{q}\right\Vert _{H^{4}\left( Q_{T}\right)
}^{2}+\right) \exp \left[ 2\lambda \left( 2A_{1}+2\right) ^{\nu _{0}}\right]
+ \\ 
+C_{1}\left( \left\Vert \widetilde{F}\right\Vert _{H^{4}\left( Q_{T}\right)
}^{2}+\left\Vert \widetilde{G}\right\Vert _{H^{4}\left( Q_{T}\right)
}^{2}\right) \exp \left[ 2\lambda \left( 2A_{1}+2\right) ^{\nu _{0}}\right] .%
\end{array}%
\right.  \label{5.23}
\end{equation}%
The second above mentioned second inequality is generated by (\ref{5.18})
and the first estimate (\ref{5.210}). This inequality is:%
\begin{equation}
\left. 
\begin{array}{c}
\dint\limits_{Q_{T}}\left( w_{t}-\Delta w\right) ^{2}\varphi _{\lambda ,\nu
_{0}}\left( x\right) dxdt\leq C_{1}\dint\limits_{Q_{T}}\left( \left( \nabla
w\right) ^{2}+w^{2}\right) \varphi _{\lambda ,\nu _{0}}\left( x\right) dxdt+
\\ 
+C_{1}\dint\limits_{Q_{T}}\left( \left( \nabla v\right) ^{2}+v^{2}\right)
\varphi _{\lambda ,\nu _{0}}\left( x\right)
dxdt+C_{1}\dint\limits_{Q_{T}}\left( \Delta v\right) ^{2}\varphi _{\lambda
,\nu _{0}}\left( x\right) dxdt+ \\ 
+C_{1}\left( \left\Vert \widetilde{p}\right\Vert _{H^{4}\left( \Omega
\right) }^{2}+\left\Vert \widetilde{q}\right\Vert _{H^{4}\left( \Omega
\right) }^{2}\right) \exp \left[ 2\lambda \left( 2A_{1}+2\right) ^{\nu _{0}}%
\right] + \\ 
+C_{1}\left( \left\Vert \widetilde{F}\right\Vert _{H^{4}\left( \Omega
\right) }^{2}+\left\Vert \widetilde{G}\right\Vert _{H^{4}\left( \Omega
\right) }^{2}\right) \exp \left[ 2\lambda \left( 2A_{1}+2\right) ^{\nu _{0}}%
\right] .%
\end{array}%
\right.  \label{5.24}
\end{equation}

It follows from the last line of (\ref{4.6}) and (\ref{5.16}) that we can
apply Carleman estimates of Theorems 4.3 and 4.2 to the left hand sides of (%
\ref{5.23}) and (\ref{5.24}) respectively. Let $\lambda _{0}>1$ be the
parameter of Theorems 4.1-4.3. Hence, using (\ref{4.6}), (\ref{5.20}) and (%
\ref{5.21}), we again obtain two estimates for all $\lambda \geq \lambda
_{0} $. The first estimate is:%
\begin{equation}
\left. 
\begin{array}{c}
\left( 1/\lambda \right) \dint\limits_{Q_{T}}\left(
v_{t}^{2}+\dsum\limits_{i,j=2}^{n}v_{ij}^{2}\right) \varphi _{\lambda ,\nu
_{0}}dxdt+\dint\limits_{Q_{T}}\left( \lambda \left( \nabla v\right)
^{2}+\lambda ^{3}v^{2}\right) \varphi _{\lambda ,\nu _{0}}dxdt\leq \\ 
\leq C_{1}\dint\limits_{Q_{T}}\left( \left( \nabla v\right)
^{2}+v^{2}\right) \varphi _{\lambda ,\nu _{0}}\left( x\right)
dxdt+C_{1}\dint\limits_{Q_{T}}w^{2}\varphi _{\lambda ,\nu _{0}}\left(
x\right) dxdt+D,%
\end{array}%
\right.  \label{5.25}
\end{equation}%
where%
\begin{equation}
\left. 
\begin{array}{c}
D= \\ 
\begin{array}{c}
=C_{1}\left( \left\Vert \widetilde{p}\right\Vert _{H^{4}\left( \Omega
\right) }^{2}+\left\Vert \widetilde{q}\right\Vert _{H^{4}\left( \Omega
\right) }^{2}\right) \exp \left[ 3\lambda \left( 2A_{1}+2\right) ^{\nu _{0}}%
\right] + \\ 
+C_{1}\left( \left\Vert \widetilde{F}\right\Vert _{H^{4}\left( \Omega
\right) }^{2}+\left\Vert \widetilde{G}\right\Vert _{H^{4}\left( \Omega
\right) }^{2}\right) \exp \left[ 3\lambda \left( 2A_{1}+2\right) ^{\nu _{0}}%
\right] + \\ 
+C_{1}\left( \left\Vert \widetilde{f}_{0}\right\Vert _{H^{2,1}\left(
S_{T}\diagdown \Gamma _{T}^{-}\right) }^{2}+\left\Vert \widetilde{f}%
_{1}\right\Vert _{H^{1,0}\left( S_{T}\diagdown \Gamma _{T}^{-}\right)
}^{2}\right) \exp \left[ 3\lambda \left( 2A_{1}+2\right) ^{\nu _{0}}\right] +
\\ 
+C_{1}\left( \left\Vert \widetilde{g}_{0}\right\Vert _{H^{2,1}\left(
S_{T}\diagdown \Gamma _{T}^{-}\right) }^{2}+\left\Vert \widetilde{g}%
_{1}\right\Vert _{H^{2,1}\left( S_{T}\diagdown \Gamma _{T}^{-}\right)
}^{2}\right) \exp \left[ 3\lambda \left( 2A_{1}+2\right) ^{\nu _{0}}\right] +
\\ 
+C_{1}\exp \left[ 3\cdot 2^{\nu _{0}}\lambda \right] \left( \left\Vert
v\right\Vert _{H^{2,1}\left( \Gamma _{T}^{-}\right) }^{2}+\left\Vert
v_{x_{1}}\right\Vert _{H^{1,0}\left( \Gamma _{T}^{-}\right) }^{2}\right) +
\\ 
+\backslash +C_{1}\exp \left[ 3\cdot 2^{\nu _{0}}\lambda \right] \left(
\left\Vert w\right\Vert _{H^{2,1}\left( \Gamma _{T}^{-}\right)
}^{2}+\left\Vert w_{x_{1}}\right\Vert _{H^{1,0}\left( \Gamma _{T}^{-}\right)
}^{2}\right) .%
\end{array}%
\end{array}%
\right.  \label{5.26}
\end{equation}%
Here we make the estimate for $D$ slightly stronger for a convenience of
further derivations. Since $C_{1}$ denotes different numbers (Remark 4.3),
then below $D$ denotes different numbers with the same expression (\ref{5.26}%
).

The second estimate is:%
\begin{equation}
\left. 
\begin{array}{c}
\left( 1/\lambda \right) \dint\limits_{Q_{T}}\left(
w_{t}^{2}+\dsum\limits_{i,j=2}^{n}w_{ij}^{2}\right) \varphi _{\lambda ,\nu
_{0}}dxdt+\dint\limits_{Q_{T}}\left( \lambda \left( \nabla w\right)
^{2}+\lambda ^{3}w^{2}\right) \varphi _{\lambda ,\nu _{0}}dxdt\leq \\ 
\leq C_{1}\dint\limits_{Q_{T}}\left( \left( \nabla w\right)
^{2}+w^{2}\right) \varphi _{\lambda ,\nu _{0}}\left( x\right) dxdt+ \\ 
+C_{1}\dint\limits_{Q_{T}}\left( \left( \nabla v\right) ^{2}+v^{2}\right)
\varphi _{\lambda ,\nu _{0}}\left( x\right)
dxdt+C_{1}\dint\limits_{Q_{T}}\left( \Delta v\right) ^{2}\varphi _{\lambda
,\nu _{0}}\left( x\right) dxdt+D.%
\end{array}%
\right.  \label{5.27}
\end{equation}%
Choose a sufficiently large $\lambda _{1}=\lambda _{1}\left( M,c,\Omega
,T\right) >\lambda _{0}$ such that for all $\lambda \geq \lambda _{1}$ and
for all functions $h\in H^{1,0}\left( Q_{T}\right) $ 
\begin{equation*}
\left. C_{1}\dint\limits_{Q_{T}}\left( \left( \nabla h\right)
^{2}+h^{2}\right) \varphi _{\lambda ,\nu _{0}}\left( x\right) dxdt\leq
\left( 1/2\right) \dint\limits_{Q_{T}}\left( \lambda \left( \nabla h\right)
^{2}+\lambda ^{3}h^{2}\right) \varphi _{\lambda ,\nu _{0}}dxdt.\right.
\end{equation*}%
The (\ref{5.25}) implies%
\begin{equation}
\left. 
\begin{array}{c}
\left( 1/\lambda \right) \dint\limits_{Q_{T}}\left(
v_{t}^{2}+\dsum\limits_{i,j=2}^{n}v_{ij}^{2}\right) \varphi _{\lambda ,\nu
_{0}}dxdt+\dint\limits_{Q_{T}}\left( \lambda \left( \nabla v\right)
^{2}+\lambda ^{3}v^{2}\right) \varphi _{\lambda ,\nu _{0}}dxdt\leq \\ 
\leq C_{1}\dint\limits_{Q_{T}}w^{2}\varphi _{\lambda ,\nu _{0}}\left(
x\right) dxdt+D,\text{ }\forall \lambda \geq \lambda _{1}.%
\end{array}%
\right.  \label{5.29}
\end{equation}%
Similarly, (\ref{5.27})\ implies%
\begin{equation}
\left. 
\begin{array}{c}
\left( 1/\lambda \right) \dint\limits_{Q_{T}}\left(
w_{t}^{2}+\dsum\limits_{i,j=2}^{n}w_{ij}^{2}\right) \varphi _{\lambda ,\nu
_{0}}dxdt+\dint\limits_{Q_{T}}\left( \lambda \left( \nabla w\right)
^{2}+\lambda ^{3}w^{2}\right) \varphi _{\lambda ,\nu _{0}}dxdt\leq \\ 
\leq C_{1}\dint\limits_{Q_{T}}\left( \left( \nabla v\right)
^{2}+v^{2}\right) \varphi _{\lambda ,\nu _{0}}\left( x\right)
dxdt+C_{1}\dint\limits_{Q_{T}}\left( \Delta v\right) ^{2}\varphi _{\lambda
,\nu _{0}}\left( x\right) dxdt+D_{2}, \\ 
\forall \lambda \geq \lambda _{1}.%
\end{array}%
\right.  \label{5.30}
\end{equation}%
Divide (\ref{5.30}) by $\lambda ^{2}$ and sum up with (\ref{5.29}). We obtain%
\begin{equation}
\left. 
\begin{array}{c}
\left( 1/\lambda \right) \dint\limits_{Q_{T}}\left(
v_{t}^{2}+\dsum\limits_{i,j=2}^{n}v_{ij}^{2}\right) \varphi _{\lambda ,\nu
_{0}}dxdt+\left( 1/\lambda ^{3}\right) \dint\limits_{Q_{T}}\left(
w_{t}^{2}+\dsum\limits_{i,j=2}^{n}w_{ij}^{2}\right) \varphi _{\lambda ,\nu
_{0}}dxdt+ \\ 
+\dint\limits_{Q_{T}}\left( \lambda \left( \nabla v\right) ^{2}+\lambda
^{3}v^{2}\right) \varphi _{\lambda ,\nu _{0}}dxdt+\dint\limits_{Q_{T}}\left[
\left( 1/\lambda \right) \left( \nabla w\right) ^{2}+\lambda w^{2}\right]
\varphi _{\lambda ,\nu _{0}}dxdt+ \\ 
+\left( C_{1}/\lambda ^{2}\right) \dint\limits_{Q_{T}}\left( \Delta v\right)
^{2}\varphi _{\lambda ,\nu _{0}}\left( x\right) dxdt+D,\text{ }\forall
\lambda \geq \lambda _{1}.%
\end{array}%
\right.  \label{5.31}
\end{equation}%
Since $\lambda _{1}$ is sufficiently large, then%
\begin{equation*}
\frac{C_{1}}{\lambda ^{2}}\dint\limits_{Q_{T}}\left( \Delta v\right)
^{2}\varphi _{\lambda ,\nu _{0}}\left( x\right) dxdt\leq \frac{1}{2\lambda }%
\dint\limits_{Q_{T}}\left(
v_{t}^{2}+\dsum\limits_{i,j=2}^{n}v_{ij}^{2}\right) \varphi _{\lambda ,\nu
_{0}}dxdt,\text{ }\forall \lambda \geq \lambda _{1}.
\end{equation*}%
Hence, (\ref{5.31}) can be rewritten as:%
\begin{equation}
\left. 
\begin{array}{c}
\left( 1/\lambda \right) \dint\limits_{Q_{T}}\left(
v_{t}^{2}+\dsum\limits_{i,j=2}^{n}v_{ij}^{2}\right) \varphi _{\lambda ,\nu
_{0}}dxdt+\left( 1/\lambda ^{3}\right) \dint\limits_{Q_{T}}\left(
w_{t}^{2}+\dsum\limits_{i,j=2}^{n}w_{ij}^{2}\right) \varphi _{\lambda ,\nu
_{0}}dxdt+ \\ 
+\dint\limits_{Q_{T}}\left( \lambda \left( \nabla v\right) ^{2}+\lambda
^{3}v^{2}\right) \varphi _{\lambda ,\nu _{0}}dxdt+\dint\limits_{Q_{T}}\left[
\left( 1/\lambda \right) \left( \nabla w\right) ^{2}+\lambda \left( w\right)
^{2}\right] \varphi _{\lambda ,\nu _{0}}dxdt+D, \\ 
\forall \lambda \geq \lambda _{1},%
\end{array}%
\right.  \label{5.33}
\end{equation}

Using (\ref{4.12}), (\ref{5.1}) and (\ref{5.26}), we obtain that, similarly
with (\ref{5.21}), the transformation procedure of sub-subsection 5.1.1
leads to the following estimate for $D$ 
\begin{equation}
\left. 
\begin{array}{c}
D\leq C_{1}\delta ^{2}\exp \left[ 3\lambda \left( 2A_{1}+2\right) ^{\nu _{0}}%
\right] + \\ 
+C_{1}\exp \left[ 3\cdot 2^{\nu _{0}}\lambda \right] \left( \left\Vert
v\right\Vert _{H^{2,1}\left( \Gamma _{T}^{-}\right) }^{2}+\left\Vert
v_{x_{1}}\right\Vert _{H^{1,0}\left( \Gamma _{T}^{-}\right) }^{2}\right) +
\\ 
+C_{1}\exp \left[ 3\cdot 2^{\nu _{0}}\lambda \right] \left( \left\Vert
w\right\Vert _{H^{2,1}\left( \Gamma _{T}^{-}\right) }^{2}+\left\Vert
w_{x_{1}}\right\Vert _{H^{1,0}\left( \Gamma _{T}^{-}\right) }^{2}\right) .%
\end{array}%
\right.  \label{5.35}
\end{equation}%
Next, using trace theorem and again similarly with (\ref{5.21}), we obtain%
\begin{equation*}
\left\Vert v\right\Vert _{H^{2,1}\left( \Gamma _{T}^{-}\right)
}^{2}+\left\Vert v_{x_{1}}\right\Vert _{H^{1,0}\left( \Gamma _{T}^{-}\right)
}^{2}+\left\Vert w\right\Vert _{H^{2,1}\left( \Gamma _{T}^{-}\right)
}^{2}+\left\Vert w_{x_{1}}\right\Vert _{H^{1,0}\left( \Gamma _{T}^{-}\right)
}^{2}\leq C_{1}.
\end{equation*}%
Hence, (\ref{5.35}) imply that $D$ can be estimated as: 
\begin{equation}
D\leq C_{1}\exp \left[ 3\lambda \left( 2A_{1}+2\right) ^{\nu _{0}}\right]
\delta ^{2}+C_{1}\exp \left[ 3\cdot 2^{\nu _{0}}\lambda \right] .
\label{5.36}
\end{equation}%
\begin{equation*}
\left\Vert v\right\Vert _{H^{2,1}\left( \Gamma _{T}^{-}\right)
}^{2}+\left\Vert v_{x_{1}}\right\Vert _{H^{1,0}\left( \Gamma _{T}^{-}\right)
}^{2}+\left\Vert w\right\Vert _{H^{2,1}\left( \Gamma _{T}^{-}\right)
}^{2}+\left\Vert w_{x_{1}}\right\Vert _{H^{1,0}\left( \Gamma _{T}^{-}\right)
}^{2}\leq C_{1}.
\end{equation*}

We now recall the domain $Q_{\gamma T}$ in (\ref{2.1}), where $\gamma \in
\left( 0,2A_{1}\right) $ is an arbitrary but fixed number. By (\ref{4.1}) $%
\varphi _{\lambda ,\nu _{0}}\left( x\right) \geq \exp \left( 2\lambda \left(
\gamma +2\right) ^{\nu _{0}}\right) $ in $Q_{\gamma T}.$ Hence, 
\begin{equation}
\exp \left( 2\lambda \left( \gamma +2\right) ^{\nu _{0}}\right) \left\Vert
f\right\Vert _{L_{2}\left( Q_{\gamma T}\right) }^{2}\leq
\dint\limits_{Q_{T}}f^{2}\left( x,t\right) \varphi _{\lambda ,\nu _{0}}dxdt,%
\text{ }\forall f\in L_{2}\left( Q_{T}\right) .  \label{5.37}
\end{equation}%
Hence, (\ref{5.33}), (\ref{5.36}) and (\ref{5.37}) imply%
\begin{equation}
\left. 
\begin{array}{c}
\left\Vert v\right\Vert _{H^{2,1}\left( Q_{\gamma T}\right) }+\left\Vert
w\right\Vert _{H^{2,1}\left( Q_{\gamma T}\right) }\leq \\ 
\leq C_{1}\left[ \exp \left( 1.5\lambda \left( 2A_{1}+2\right) ^{\nu
_{0}}\right) \delta +\exp \left( -1.5\lambda \left( \left( \gamma +2\right)
^{\nu _{0}}-2^{\nu _{0}}\right) \right) \right] , \\ 
\forall \lambda \geq \lambda _{1}.%
\end{array}%
\right.  \label{5.38}
\end{equation}

Choose now $\lambda =\lambda \left( \delta \right) $ such that (\ref{5.38})%
\begin{equation*}
\exp \left( 1.5\lambda \left( 2A_{1}+2\right) ^{\nu _{0}}\right) \delta
=\exp \left( -1.5\lambda \left( \left( \gamma +2\right) ^{\nu _{0}}-2^{\nu
_{0}}\right) \right) .
\end{equation*}%
Hence,%
\begin{equation}
1.5\left[ \left( 2A_{1}+2\right) ^{\nu _{0}}+\left( \gamma +2\right) ^{\nu
_{0}}-2^{\nu _{0}}\right] \lambda =\ln \left( \delta ^{-1}\right) .
\label{5.39}
\end{equation}%
Recall that by Theorem 4.1, $\nu _{0}=\nu _{0}\left( A_{1}\right) .$ Also,
recall that $\lambda _{1}=\lambda _{1}\left( M,c,\gamma ,\Omega ,T\right) .$
Hence, by (\ref{5.39}) 
\begin{equation}
\left. 
\begin{array}{c}
\lambda \left( \delta \right) =\ln \left( \delta ^{-1/d}\right) , \\ 
d=1.5\left[ \left( \gamma +2\right) ^{\nu _{0}}-2^{\nu _{0}}+\left(
2A_{1}+2\right) ^{\nu _{0}}\right] , \\ 
\forall \delta \in \left( 0,\delta _{0}\right) ,\text{ }\delta _{0}=\delta
_{0}\left( M,c,\Omega ,T\right) :\ln \left( \delta _{0}^{-1/d}\right)
>\lambda _{1}.%
\end{array}%
\right.  \label{5.40}
\end{equation}%
Hence, by Remark 4.3, we should now replace $C_{1}=C_{1}\left( M,c,\Omega
,T\right) >0$ with $B=B\left( M,c,\gamma ,\Omega ,T\right) >0$\emph{.}
Consider the number $\alpha \in \left( 0,1\right) ,$%
\begin{equation}
\alpha \left( M,c,\gamma ,\Omega ,T\right) =\frac{1.5\left( \left( \gamma
+2\right) ^{\nu _{0}}-2^{\nu _{0}}\right) }{d}=\frac{\left( \gamma +2\right)
^{\nu _{0}}-2^{\nu _{0}}}{\left( \gamma +2\right) ^{\nu _{0}}-2^{\nu
_{0}}+\left( 2A_{1}+2\right) ^{\nu _{0}}}.  \label{5.41}
\end{equation}%
It follows from (\ref{5.38})-(\ref{5.41}) that the following H\"{o}lder
stability estimate for functions $v$ and $w$ is valid:%
\begin{equation}
\left\Vert v\right\Vert _{H^{2,1}\left( Q_{\gamma T}\right) }+\left\Vert
w\right\Vert _{H^{2,1}\left( Q_{\gamma T}\right) }\leq B\delta ^{\alpha },%
\text{ }\forall \delta \in \left( 0,\delta _{0}\right) .  \label{5.42}
\end{equation}%
It follows from (\ref{5.16}), (\ref{5.42}) and trace theorem that estimate (%
\ref{4.15}) is valid, which is the second target estimate of this theorem.

To prove target estimates (\ref{4.14}) of this theorem, we recall again the
transformation procedure of sub-subsection 5.1.1. Using (\ref{4.12}), (\ref%
{5.1}), (\ref{5.2}), (\ref{5.6})-(\ref{5.15}) and triangle inequality, we
obtain%
\begin{equation}
\left. 
\begin{array}{c}
\left\Vert v\right\Vert _{H^{2,1}\left( Q_{\gamma T}\right) }+\left\Vert
w\right\Vert _{H^{2,1}\left( Q_{\gamma T}\right) }\geq \left\Vert \partial
_{t}^{j}u_{1}-\partial _{t}^{j}u_{2}\right\Vert _{H^{2,1}\left( Q_{\gamma
T}\right) }+ \\ 
+\left\Vert \partial _{t}^{j}m_{1}-\partial _{t}^{j}m_{2}\right\Vert
_{H^{2,1}\left( Q_{\gamma T}\right) }- \\ 
-\left( \left\Vert \widetilde{p}\right\Vert _{H^{4}\left( \Omega \right)
}+\left\Vert \widetilde{q}\right\Vert _{H^{4}\left( \Omega \right)
}+\left\Vert \widetilde{F}\right\Vert _{H^{4}\left( \Omega \right)
}+\left\Vert \widetilde{G}\right\Vert _{H^{4}\left( \Omega \right) }\right)
\geq \\ 
\geq \left\Vert \partial _{t}^{j}u_{1}-\partial _{t}^{j}u_{2}\right\Vert
_{H^{2,1}\left( Q_{\gamma T}\right) }+\left\Vert \partial
_{t}^{j}m_{1}-\partial _{t}^{j}m_{2}\right\Vert _{H^{2,1}\left( Q_{\gamma
T}\right) }-C_{1}\delta ,\text{ }j=0,1.%
\end{array}%
\right.  \label{5.240}
\end{equation}%
Comparing this with (\ref{5.42}) and using $\delta <\delta ^{\alpha },$ $%
\forall \delta \in \left( 0,\delta _{0}\right) ,$ we obtain (\ref{4.14}).

To prove uniqueness, we set $\delta =0.$ Then (\ref{4.14}) and (\ref{4.15})
imply that 
\begin{equation}
u_{1}\left( x,t\right) =u_{2}\left( x,t\right) =m_{1}\left( x,t\right)
=m_{2}\left( x,t\right) =0\text{ in }Q_{\gamma T}\text{ and }b_{1}\left(
x\right) =b_{2}\left( x\right) \text{ in }\Omega _{\gamma }.  \label{5.43}
\end{equation}%
Since $\gamma \in \left( 0,2A_{1}\right) $ is an arbitrary number, then,
setting $\gamma \rightarrow 0,$ we obtain that (\ref{5.43}) holds for $Q_{T}$
and $\Omega .$ $\square $

\subsection{Proof of Theorem 4.5}

\QTP{Body Math}
\label{sec:5.2}

The proof of this theorem can be carried out as an insignificant
modification of the proof of Theorem 4.4. Indeed, since by (\ref{2.6}) the
lateral Cauchy data are known now at the entire boundary $S_{T},$ then we
should not separate $S_{T}\diagdown \Gamma _{T}^{-}$ from $\Gamma _{T}^{-}$
in the above proof. In particular, first and second lines in the Carleman
estimate (\ref{4.6}) should be replaced with:%
\begin{equation}
C\exp \left[ 3\lambda \left( 2A_{1}+2\right) ^{\nu _{0}}\right] \left(
\left\Vert u\right\Vert _{H^{2,1}\left( S_{T}\right) }^{2}+\left\Vert
\partial _{n}u\right\Vert _{H^{1,0}\left( S_{T}{}^{-}\right) }^{2}\right) .
\label{5.44}
\end{equation}%
Hence, $D$ in (\ref{5.26}) should be replaced with:%
\begin{equation}
\left. 
\begin{array}{c}
D= \\ 
\begin{array}{c}
=C_{1}\left( \left\Vert \widetilde{p}\right\Vert _{H^{4}\left( \Omega
\right) }^{2}+\left\Vert \widetilde{q}\right\Vert _{H^{4}\left( \Omega
\right) }^{2}+\right) \exp \left[ 3\lambda \left( 2A_{1}+2\right) ^{\nu _{0}}%
\right] + \\ 
+C_{1}\left( \left\Vert \widetilde{F}\right\Vert _{H^{4}\left( \Omega
\right) }^{2}+\left\Vert \widetilde{G}\right\Vert _{H^{4}\left( \Omega
\right) }^{2}\right) \exp \left[ 3\lambda \left( 2A_{1}+2\right) ^{\nu _{0}}%
\right] + \\ 
+C_{1}\left( \left\Vert \widetilde{f}_{0}\right\Vert _{H^{2,1}\left(
S_{T}\right) }^{2}+\left\Vert \widetilde{f}_{1}\right\Vert _{H^{1,0}\left(
S_{T}\right) }^{2}\right) \exp \left[ 3\lambda \left( 2A_{1}+2\right) ^{\nu
_{0}}\right] + \\ 
+C_{1}\left( \left\Vert \widetilde{g}_{0}\right\Vert _{H^{2,1}\left(
S_{T}{}^{-}\right) }^{2}+\left\Vert \widetilde{g}_{1}\right\Vert
_{H^{2,1}\left( S_{T}{}^{-}\right) }^{2}\right) \exp \left[ 3\lambda \left(
2A_{1}+2\right) ^{\nu _{0}}\right] .%
\end{array}%
\end{array}%
\right.  \label{5.45}
\end{equation}%
Using (\ref{5.33}), (\ref{5.44}) and (\ref{5.45}), we obtain%
\begin{equation}
\left. 
\begin{array}{c}
\left( 1/\lambda \right) \dint\limits_{Q_{T}}\left(
v_{t}^{2}+\dsum\limits_{i,j=2}^{n}v_{ij}^{2}\right) \varphi _{\lambda ,\nu
_{0}}dxdt+\left( 1/\lambda ^{3}\right) \dint\limits_{Q_{T}}\left(
w_{t}^{2}+\dsum\limits_{i,j=2}^{n}w_{ij}^{2}\right) \varphi _{\lambda ,\nu
_{0}}dxdt+ \\ 
+\dint\limits_{Q_{T}}\left( \lambda \left( \nabla v\right) ^{2}+\lambda
^{3}v^{2}\right) \varphi _{\lambda ,\nu _{0}}dxdt+\dint\limits_{Q_{T}}\left[
\left( 1/\lambda \right) \left( \nabla w\right) ^{2}+\lambda \left( w\right)
^{2}\right] \varphi _{\lambda ,\nu _{0}}dxdt+ \\ 
+C_{1}\left( \left\Vert \widetilde{p}\right\Vert _{H^{4}\left( \Omega
\right) }^{2}+\left\Vert \widetilde{q}\right\Vert _{H^{4}\left( \Omega
\right) }^{2}\right) \exp \left[ 3\lambda \left( 2A_{1}+2\right) ^{\nu _{0}}%
\right] + \\ 
+C_{1}\left( \left\Vert \widetilde{F}\right\Vert _{H^{4}\left( \Omega
\right) }^{2}+\left\Vert \widetilde{G}\right\Vert _{H^{4}\left( \Omega
\right) }^{2}\right) \exp \left[ 3\lambda \left( 2A_{1}+2\right) ^{\nu _{0}}%
\right] + \\ 
+C_{1}\left( \left\Vert \widetilde{f}_{0}\right\Vert _{H^{2,1}\left(
S_{T}\right) }^{2}+\left\Vert \widetilde{f}_{1}\right\Vert _{H^{1,0}\left(
S_{T}\right) }^{2}\right) \exp \left[ 3\lambda \left( 2A_{1}+2\right) ^{\nu
_{0}}\right] + \\ 
++C_{1}\left( \left\Vert \widetilde{g}_{0}\right\Vert _{H^{2,1}\left(
S_{T}{}^{-}\right) }^{2}+\left\Vert \widetilde{g}_{1}\right\Vert
_{H^{2,1}\left( S_{T}{}^{-}\right) }^{2}\right) \exp \left[ 3\lambda \left(
2A_{1}+2\right) ^{\nu _{0}}\right] , \\ 
\forall \lambda \geq \lambda _{1},%
\end{array}%
\right.  \label{5.46}
\end{equation}%
Since by (\ref{4.2}) $\varphi _{\lambda ,\nu _{0}}\left( x\right) \geq \exp
\left( 2\lambda \cdot 2^{\nu }\right) ,$ then we set in (\ref{5.46}) $%
\lambda =\lambda _{1}$ divide by $\exp \left( 2\lambda \cdot 2^{\nu }\right) 
$ and then obtain similarly with (\ref{5.37}) and (\ref{5.38}):%
\begin{equation*}
\left. 
\begin{array}{c}
\left\Vert v\right\Vert _{H^{2,1}\left( Q_{T}\right) }+\left\Vert
w\right\Vert _{H^{2,1}\left( Q_{T}\right) }\leq \\ 
\leq C_{1}\left( \left\Vert \widetilde{p}\right\Vert _{H^{4}\left( \Omega
\right) }+\left\Vert \widetilde{q}\right\Vert _{H^{4}\left( \Omega \right)
}+\left\Vert \widetilde{F}\right\Vert _{H^{4}\left( \Omega \right)
}+\left\Vert \widetilde{G}\right\Vert _{H^{4}\left( \Omega \right) }\right) +
\\ 
+C_{1}\left( \left\Vert \widetilde{f}_{0}\right\Vert _{H^{2,1}\left(
S_{T}\right) }+\left\Vert \widetilde{f}_{1}\right\Vert _{H^{1,0}\left(
S_{T}\right) }+\left\Vert \widetilde{g}_{0}\right\Vert _{H^{2,1}\left(
S_{T}\right) }+\left\Vert \widetilde{g}_{1}\right\Vert _{H^{1,0}\left(
S_{T}\right) }\right) .%
\end{array}%
\right.
\end{equation*}%
The rest of the proof is similar with the proof of subsection 4.1, starting
from (\ref{5.37}). $\square $

\section{Appendix: Proofs of Theorems 4.1 and 4.2}

\label{sec:6}

\subsection{Proof of Theorem 4.1}

\label{sec:6.1}

Recall that in this theorem $u\in C^{4,2}\left( \overline{Q}_{T}\right) .$
It is convenient not to fix $\nu $ in the major part of this proof. Rather,
we assume that $\nu \geq \nu _{0}>1$ and set $\nu =\nu _{0}\left(
A_{1}\right) $ only when being close to the end of the proof. The constant $%
C $ is independent on $\nu .$ In this proof $\lambda \geq \lambda _{0}\left(
A_{1}\right) $ and both parameters $\lambda _{0},\nu _{0}$ are sufficiently
large. Furthermore since $\nu =\nu _{0}\left( A_{1}\right) $ in the end,
then we assume that 
\begin{equation}
\lambda >>\nu .  \label{7.1}
\end{equation}%
Using (\ref{4.1}), change variables%
\begin{equation}
v=ue^{\lambda \psi ^{\nu }}\rightarrow u=ve^{-\lambda \psi ^{\nu }}.
\label{7.2}
\end{equation}%
Hence,%
\begin{equation}
\left. 
\begin{array}{c}
u_{t}=v_{t}e^{-\lambda \psi ^{\nu }},\text{ }u_{x_{1}}=\left(
v_{x_{1}}-\lambda \nu \psi ^{\nu -1}v\right) e^{-\lambda \psi ^{\nu }}, \\ 
u_{x_{1}x_{1}}=\left\{ v_{x_{1}x_{1}}-2\lambda \nu \psi ^{\nu
-1}v_{x_{1}}+\lambda ^{2}\nu ^{2}\psi ^{2\nu -2}\left[ 1-2\psi ^{-\nu
}\left( \nu -1\right) /\left( \lambda \nu \right) \right] v\right\}
e^{-\lambda \psi ^{\nu }}, \\ 
u_{x_{i}x_{i}}=v_{x_{i}x_{i}}e^{-\lambda \psi ^{\nu }},\text{ }i,j=2,...,n.%
\end{array}%
\right.  \label{7.3}
\end{equation}%
By (\ref{7.2}) and (\ref{7.3})%
\begin{equation}
\left. 
\begin{array}{c}
\left( u_{t}-\Delta u\right) ^{2}\varphi _{\lambda ,\nu }\psi ^{-\nu +1}= \\ 
=\left[ 
\begin{array}{c}
v_{t}-\left( v_{x_{1}x_{1}}+\dsum\limits_{i=2}^{n}v_{x_{i}x_{i}}\right)
+2\lambda \nu \psi ^{\nu -1}v_{x_{1}}- \\ 
-\lambda ^{2}\nu ^{2}\psi ^{2\nu -2}\left[ 1-2\psi ^{-\nu }\left( \nu
-1\right) /\left( \lambda \nu \right) \right] v%
\end{array}%
\right] ^{2}\psi ^{-\nu +1}.%
\end{array}%
\right.  \label{7.4}
\end{equation}%
Denote%
\begin{equation}
\left. 
\begin{array}{c}
z_{1}=v_{t},\text{ }z_{2}=-v_{x_{1}x_{1}}-\dsum%
\limits_{i,j=2}^{n}v_{x_{i}x_{i}}, \\ 
z_{3}=2\lambda \nu \psi ^{\nu -1}v_{x_{1}},\text{ } \\ 
z_{4}=-\lambda ^{2}\nu ^{2}\psi ^{2\nu -2}\left[ 1-2\psi ^{-\nu }\left( \nu
-1\right) /\left( \lambda \nu \right) \right] v.%
\end{array}%
\right.  \label{7.5}
\end{equation}%
By (\ref{7.4}) and (\ref{7.5})%
\begin{equation}
\left. 
\begin{array}{c}
\left( u_{t}-\Delta u\right) ^{2}\varphi _{\lambda ,\nu }\psi ^{-\nu +1}= 
\left[ \left( z_{1}+z_{3}\right) +z_{2}+z_{4}\right] ^{2}\psi ^{-\nu +1}\geq
\\ 
\geq \left( z_{1}+z_{3}\right) ^{2}\psi ^{-\nu +1}+2z_{1}z_{2}\psi ^{-\nu
+1}+2z_{1}z_{4}\psi ^{-\nu +1}+2z_{2}z_{3}\psi ^{-\nu +1}+2z_{3}z_{4}\psi
^{-\nu +1}.%
\end{array}%
\right.  \label{7.6}
\end{equation}

\subsubsection{Step 1. Estimate from the below the term $2z_{1}z_{2}\protect%
\psi ^{-\protect\nu +1}$ in (\protect\ref{7.6})}

We have: 
\begin{equation*}
\left. 
\begin{array}{c}
2z_{1}z_{2}\psi ^{-\nu +1}=-2v_{t}v_{x_{1}x_{1}}\psi ^{-\nu
+1}-2\dsum\limits_{i,j=2}^{n}v_{t}v_{x_{i}x_{i}}\psi ^{-\nu +1}= \\ 
=\left( -2v_{t}v_{x_{1}}\psi ^{-\nu +1}\right)
_{x_{1}}+2v_{tx_{1}}v_{x_{1}}\psi ^{-\nu +1}-2\left( \nu -1\right) \psi
^{-\nu }v_{t}v_{x_{1}}+ \\ 
+\dsum\limits_{i=2}^{n}\left( -2v_{t}v_{x_{i}}\psi ^{-\nu +1}\right)
_{x_{i}}+\dsum\limits_{i=2}^{n}2v_{tx_{i}}v_{x_{i}}\psi ^{-\nu +1}= \\ 
=-2\left( \nu -1\right) \psi ^{-\nu }z_{1}v_{x_{1}}+ \\ 
+\left( -2v_{t}v_{x_{1}}\psi ^{-\nu +1}\right)
_{x_{1}}+\dsum\limits_{i=2}^{n}\left( -2v_{t}v_{x_{i}}\psi ^{-\nu +1}\right)
_{x_{i}}+ \\ 
+\partial _{t}\left( v_{x_{1}}^{2}\psi ^{-\nu
+1}+\dsum\limits_{i,j=2}^{n}v_{x_{i}}^{2}\psi ^{-\nu +1}\right) .%
\end{array}%
\right.
\end{equation*}%
Thus, 
\begin{equation}
2z_{1}z_{2}\psi ^{-\nu +1}=-2\left( \nu -1\right) \psi ^{-\nu
}z_{1}v_{x_{1}}+\partial _{t}V_{1}+\func{div}U_{1},  \label{7.7}
\end{equation}%
where%
\begin{equation}
\partial _{t}V_{1}=\partial _{t}\left[ \left( \left( u_{x_{1}}+\lambda \nu
\psi ^{\nu -1}u\right) ^{2}+\dsum\limits_{i,j=2}^{n}u_{x_{i}}^{2}\right)
\varphi _{\lambda ,\nu }\psi ^{-\nu +1}\right] ,  \label{7.8}
\end{equation}%
\begin{equation}
\func{div}U_{1}=\left( -2u_{t}\left( u_{x_{1}}+\lambda \nu \psi ^{\nu
-1}u\right) \varphi _{\lambda ,\nu }\psi ^{-\nu +1}\right)
_{x_{1}}+\dsum\limits_{i=2}^{n}\left( -2u_{t}u_{x_{i}}\varphi _{\lambda ,\nu
}\psi ^{-\nu +1}\right) _{x_{i}}.  \label{7.9}
\end{equation}

\subsubsection{Step 2. Using (\protect\ref{7.5}) and (\protect\ref{7.7})-(%
\protect\ref{7.9}), estimate from the below the term $\left(
z_{1}+z_{3}\right) ^{2}\protect\psi ^{-\protect\nu +1}+2z_{1}z_{2}\protect%
\psi ^{-\protect\nu +1}$ in (\protect\ref{7.6})}

We have:%
\begin{equation}
\left. 
\begin{array}{c}
\left( z_{1}+z_{3}\right) ^{2}\psi ^{-\nu +1}+2z_{1}z_{2}\psi ^{-\nu +1}= \\ 
=z_{1}^{2}\psi ^{-\nu +1}+z_{3}^{2}\psi ^{-\nu +1}+2z_{1}z_{3}\psi ^{-\nu
+1}-2\left( \nu -1\right) \psi ^{-\nu }z_{1}v_{x_{1}}- \\ 
+\partial _{t}V_{1}+\func{div}U_{1}.%
\end{array}%
\right.  \label{7.10}
\end{equation}%
By (\ref{7.5})%
\begin{equation*}
-2\left( \nu -1\right) \psi ^{-\nu }z_{1}v_{x_{1}}=-\frac{\left( \nu
-1\right) }{\lambda \nu }\psi ^{-2\nu +1}z_{1}z_{3}.
\end{equation*}%
Hence, (\ref{7.10}) becomes%
\begin{equation}
\left. 
\begin{array}{c}
\left( z_{1}+z_{3}\right) ^{2}\psi ^{-\nu +1}+2z_{1}z_{2}\psi ^{-\nu +1}= \\ 
=\left[ z_{1}^{2}+2z_{1}z_{3}\left( 1-\left( \nu -1\right) /\left( \lambda
\nu \right) \psi ^{-\nu }\right) +z_{3}^{2}\right] \psi ^{-\nu +1}+ \\ 
+\partial _{t}V_{1}+\func{div}U_{1}.%
\end{array}%
\right.  \label{7.11}
\end{equation}%
Since $\left( 1-\left( \nu -1\right) /\left( \lambda \nu \right) \psi ^{-\nu
}\right) <1$ for sufficiently large $\lambda ,$ then 
\begin{equation*}
\left[ z_{1}^{2}+2z_{1}z_{3}\left( 1-\left( \nu -1\right) /\left( \lambda
\nu \right) \psi ^{-\nu }\right) +z_{3}^{2}\right] \psi ^{-\nu +1}\geq 0
\end{equation*}%
as a quadratic polynomial with respect to $z_{1},z_{3}.$ Hence, (\ref{7.11})
implies 
\begin{equation}
\left( z_{1}+z_{3}\right) ^{2}\psi ^{-\nu +1}+2z_{1}z_{2}\psi ^{-\nu +1}\geq
\partial _{t}V_{1}+\func{div}U_{1}.  \label{7.12}
\end{equation}

\subsubsection{Step 3. Using in (\protect\ref{7.5}), evaluate the term $%
2z_{2}z_{3}\protect\psi ^{-\protect\nu +1}$ in (\protect\ref{7.6})}

We have: 
\begin{equation*}
\left. 
\begin{array}{c}
2z_{2}z_{3}\psi ^{-\nu +1}=-4\lambda \nu v_{x_{1}}\left(
v_{x_{1}x_{1}}+\dsum\limits_{i=2}^{n}v_{x_{i}x_{i}}\right) = \\ 
=\left( -2\lambda \nu v_{x_{1}}^{2}\right)
_{x_{1}}+\dsum\limits_{i=2}^{n}\left( -4\lambda \nu
v_{x_{1}}v_{x_{i}}\right) _{x_{i}}+\dsum\limits_{i=2}^{n}\left( 4\lambda \nu
v_{x_{i}x_{1}}v_{x_{i}}\right) = \\ 
=\left( -2\lambda \nu v_{x_{1}}^{2}\right)
_{x_{1}}+\dsum\limits_{i=2}^{n}\left( -4\lambda \nu
v_{x_{1}}v_{x_{i}}\right) _{x_{i}}+\left( 2\lambda \nu
\dsum\limits_{i=2}^{n}v_{x_{i}}^{2}\right) _{x_{1}}.%
\end{array}%
\right.
\end{equation*}%
Hence,%
\begin{equation}
\left. 
\begin{array}{c}
2z_{2}z_{3}\psi ^{-\nu +1}=\func{div}U_{2}, \\ 
\func{div}U_{2}=\left( -2\lambda \nu \left( u_{x_{1}}+\lambda \nu \psi ^{\nu
-1}u\right) ^{2}\varphi _{\lambda ,\nu }+2\lambda \nu
\dsum\limits_{i=2}^{n}u_{x_{i}}^{2}\varphi _{\lambda ,\nu }\right) _{x_{1}}+
\\ 
+\dsum\limits_{i=2}^{n}\left( -4\lambda \nu \left( u_{x_{1}}+\lambda \nu
\psi ^{\nu -1}u\right) u_{x_{i}}\varphi _{\lambda ,\nu }\right) _{x_{i}}.%
\end{array}%
\right.  \label{7.13}
\end{equation}

\subsubsection{Step 4. Using in (\protect\ref{7.5}), estimate from the below
the term $2z_{1}z_{4}\protect\psi ^{-\protect\nu +1}+2z_{3}z_{4}\protect\psi %
^{-\protect\nu +1}$ in (\protect\ref{7.6})}

We have:%
\begin{equation*}
\left. 
\begin{array}{c}
2z_{1}z_{4}\psi ^{-\nu +1}+2z_{3}z_{4}\psi ^{-\nu +1}= \\ 
=-2\lambda ^{2}\nu ^{2}\psi ^{\nu -1}\left( 1-2\psi ^{-\nu }\left( \nu
-1\right) /\left( \lambda \nu \right) \right) v_{t}v- \\ 
-4\lambda ^{3}\nu ^{3}\psi ^{2\nu -2}\left( 1-2\psi ^{-\nu }\left( \nu
-1\right) /\left( \lambda \nu \right) \right) v_{x_{1}}v= \\ 
\geq \partial _{t}\left( -\lambda ^{2}\nu ^{2}\psi ^{\nu -1}\left( 1-2\psi
^{-\nu }\left( \nu -1\right) /\left( \lambda \nu \right) \right)
v^{2}\right) + \\ 
+\partial _{x_{1}}\left( -2\lambda ^{3}\nu ^{3}\psi ^{2\nu -2}\left( 1-2\psi
^{-\nu }\left( \nu -1\right) /\left( \lambda \nu \right) \right)
v^{2}\right) + \\ 
+2\lambda ^{3}\nu ^{4}\psi ^{2\nu -3}\left( 1-2\psi ^{-\nu }\left( \nu
-1\right) /\left( \lambda \nu \right) \right) v^{2}.%
\end{array}%
\right.
\end{equation*}%
We have used here the inequality $\left( \nu -1\right) >\nu /2,$ which is
valid since $\nu >2.$ Thus,%
\begin{equation}
\left. 
\begin{array}{c}
2z_{1}z_{4}\psi ^{-\nu +1}+2z_{3}z_{4}\psi ^{-\nu +1}\geq \\ 
\geq 2\lambda ^{3}\nu ^{4}\psi ^{2\nu -3}\left[ 1-2\psi ^{-\nu }\left( \nu
-1\right) /\left( \lambda \nu \right) \right] v^{2}+ \\ 
+\partial _{t}V_{2}+\func{div}U_{3},%
\end{array}%
\right.  \label{7.14}
\end{equation}%
\begin{equation}
\partial _{t}V_{2}=\partial _{t}\left( -\lambda ^{2}\nu ^{2}\psi ^{\nu
-1}\left( 1-2\psi ^{-\nu }\left( \nu -1\right) /\left( \lambda \nu \right)
\right) u^{2}\varphi _{\lambda ,\nu }\right) ,  \label{7.15}
\end{equation}%
\begin{equation}
\func{div}U_{3}=\left( -2\lambda ^{3}\nu ^{3}\psi ^{2\nu -2}\left( 1-2\psi
^{-\nu }\left( \nu -1\right) /\left( \lambda \nu \right) \right)
u^{2}\varphi _{\lambda ,\nu }\right) _{x_{1}}.  \label{7.16}
\end{equation}

\subsubsection{Step 5. Sum up (\protect\ref{7.8}), (\protect\ref{7.9}) and (%
\protect\ref{7.12})-(\protect\ref{7.16})}

Then, comparing with (\ref{7.4}) and (\ref{7.5}), we obtain 
\begin{equation}
\left. 
\begin{array}{c}
\left( u_{t}-\Delta u\right) ^{2}\varphi _{\lambda ,\nu }\psi ^{-\nu +1}\geq
\\ 
\geq 2\lambda ^{3}\nu ^{4}\psi ^{2\nu -3}\left[ 1-2\psi ^{-\nu }\left( \nu
-1\right) /\left( \lambda \nu \right) \right] u^{2}\varphi _{\lambda ,\nu }+
\\ 
+\partial _{t}V_{3}+\func{div}U_{4},%
\end{array}%
\right.  \label{7.17}
\end{equation}%
\begin{equation}
\left. 
\begin{array}{c}
\partial _{t}V_{3}=\partial _{t}V_{1}+\partial _{t}V_{2}= \\ 
=\partial _{t}\left[ \left( \left( u_{x_{1}}+\lambda \nu \psi ^{\nu
-1}u\right) ^{2}+\dsum\limits_{i,j=2}^{n}u_{x_{i}}^{2}\right) \varphi
_{\lambda ,\nu }\psi ^{-\nu +1}\right] + \\ 
+\partial _{t}\left( -\lambda ^{2}\nu ^{2}\psi ^{\nu -1}\left( 1-2\psi
^{-\nu }\left( \nu -1\right) /\left( \lambda \nu \right) \right)
u^{2}\varphi _{\lambda ,\nu }\right) ,%
\end{array}%
\right.  \label{7.18}
\end{equation}%
\begin{equation}
\left. 
\begin{array}{c}
\func{div}U_{4}=\func{div}U_{1}+\func{div}U_{2}+\func{div}U_{3}= \\ 
=\left( -2u_{t}\left( u_{x_{1}}+\lambda \nu \psi ^{\nu -1}u\right) \varphi
_{\lambda ,\nu }\psi ^{-\nu +1}\right) _{x_{1}} \\ 
+\left( -2\lambda \nu \left( u_{x_{1}}+\lambda \nu \psi ^{\nu -1}u\right)
^{2}\varphi _{\lambda ,\nu }+2\lambda \nu
\dsum\limits_{i=2}^{n}u_{x_{i}}^{2}\varphi _{\lambda ,\nu }\right) _{x_{1}}+
\\ 
+\left( -2\lambda ^{3}\nu ^{3}\psi ^{2\nu -2}\left( 1-2\psi ^{-\nu }\left(
\nu -1\right) /\left( \lambda \nu \right) \right) u^{2}\varphi _{\lambda
,\nu }\right) _{x_{1}} \\ 
+\dsum\limits_{i=2}^{n}\left( -4\lambda \nu \left( u_{x_{1}}+\lambda \nu
\psi ^{\nu -1}u\right) u_{x_{i}}\varphi _{\lambda ,\nu
}-2u_{t}u_{x_{i}}\varphi _{\lambda ,\nu }\psi ^{-\nu +1}\right) _{x_{i}}.%
\end{array}%
\right.  \label{7.19}
\end{equation}

\subsubsection{Step 6. Evaluate $\left( u_{t}-\Delta u\right) u\protect%
\varphi _{\protect\lambda ,\protect\nu }$}

We have:%
\begin{equation}
\left. 
\begin{array}{c}
\left( u_{t}-\Delta u\right) u\varphi _{\lambda ,\nu }=\partial _{t}\left(
\left( 1/2\right) u^{2}\varphi _{\lambda ,\nu }\right)
-u_{x_{1}x_{1}}u\varphi _{\lambda ,\nu
}-\dsum\limits_{i=2}^{n}u_{x_{i}x_{i}}u\varphi _{\lambda ,\nu }= \\ 
=\left( -u_{x_{1}}u\varphi _{\lambda ,\nu }\right)
_{x_{1}}+u_{x_{1}}^{2}\varphi _{\lambda ,\nu }+2\lambda \nu \psi ^{\nu
-1}u_{x_{1}}u\varphi _{\lambda ,\nu }+ \\ 
+\dsum\limits_{i=2}^{n}\left( -u_{x_{i}}u\varphi _{\lambda ,\nu }\right)
_{x_{i}}+\dsum\limits_{i=2}^{n}u_{x_{i}}^{2}\varphi _{\lambda ,\nu
}+\partial _{t}\left( \left( 1/2\right) u^{2}\varphi _{\lambda ,\nu }\right)
= \\ 
=\left( \nabla u\right) ^{2}\varphi _{\lambda ,\nu }+\left(
-u_{x_{1}}u\varphi _{\lambda ,\nu }+\lambda \nu \psi ^{\nu -1}u^{2}\varphi
_{\lambda ,\nu }\right) _{x_{1}}- \\ 
-4\lambda ^{2}\nu ^{2}\psi ^{2\nu -2}u^{2}\varphi _{\lambda ,\nu }-\lambda
\nu \left( \nu -1\right) \psi ^{\nu -2}u^{2}\varphi _{\lambda ,\nu }+ \\ 
+\dsum\limits_{i=2}^{n}\left( -u_{x_{i}}u\varphi _{\lambda ,\nu }\right)
_{x_{i}}+\partial _{t}\left( \left( 1/2\right) u^{2}\varphi _{\lambda ,\nu
}\right) .%
\end{array}%
\right.  \label{7.20}
\end{equation}

\bigskip Thus, since by (\ref{4.2}) $\lambda ^{2}\nu ^{2}\psi ^{2\nu
-2}>>\lambda \nu \left( \nu -1\right) \psi ^{\nu -2}$ for $\nu >2$ and
sufficiently large $\lambda >1$, then (\ref{7.20}) implies that for these
values of $\nu $ and $\lambda $%
\begin{equation}
\left( u_{t}-\Delta u\right) u\varphi _{\lambda ,\nu }\geq \left( \nabla
u\right) ^{2}\varphi _{\lambda ,\nu }-C\lambda ^{2}\nu ^{2}\psi ^{2\nu
-2}u^{2}\varphi _{\lambda ,\nu }+\partial _{t}V_{4}+\func{div}U_{5},
\label{7.21}
\end{equation}%
\begin{equation}
\partial _{t}V_{4}=\partial _{t}\left( \frac{1}{2}u^{2}\varphi _{\lambda
,\nu }\right) ,  \label{7.22}
\end{equation}%
\begin{equation}
\left. \func{div}U_{5}=\left( -u_{x_{1}}u\varphi _{\lambda ,\nu }+\lambda
\nu \psi ^{\nu -1}u^{2}\varphi _{\lambda ,\nu }\right)
_{x_{1}}+\dsum\limits_{i=2}^{n}\left( -u_{x_{i}}u\varphi _{\lambda ,\nu
}\right) _{x_{i}}.\right.  \label{7.23}
\end{equation}

\subsubsection{Step 7. Multiply (\protect\ref{7.21})-(\protect\ref{7.23}) by 
$\protect\lambda $ and sum up with (\protect\ref{7.17})-(\protect\ref{7.19})}

Use the inequality%
\begin{equation*}
2\lambda ^{3}\nu ^{4}\psi ^{2\nu -3}\left[ 1-2\psi ^{-\nu }\left( \nu
-1\right) /\left( \lambda \nu \right) \right] >\lambda ^{3}\nu ^{4}\psi
^{2\nu -3}>\lambda ^{3}\nu ^{2}\psi ^{2\nu -2},\forall \nu \geq \nu _{0},
\end{equation*}%
with a sufficiently large $\nu _{0}=\nu _{0}\left( A_{1}\right) >2$. We
obtain for all sufficiently large $\lambda \geq \lambda _{0}=\lambda
_{0}\left( A_{1}\right) >1$ and $\nu \geq \nu _{0}:$%
\begin{equation}
\left. 
\begin{array}{c}
\lambda \left( u_{t}-\Delta u\right) u\varphi _{\lambda ,\nu }+\left(
u_{t}-\Delta u\right) ^{2}\varphi _{\lambda ,\nu }\psi ^{-\nu +1}\geq \\ 
+\lambda \left( \nabla u\right) ^{2}\varphi _{\lambda ,\nu }+C\lambda
^{3}\nu ^{4}\psi ^{2\nu -2}u^{2}+\partial _{t}V_{5}+\func{div}U_{6},%
\end{array}%
\right.  \label{7.25}
\end{equation}%
\begin{equation}
\left. 
\begin{array}{c}
\partial _{t}V_{5}=\partial _{t}V_{3}+\partial _{t}\left( \lambda
V_{4}\right) = \\ 
=\partial _{t}\left[ \left( \left( u_{x_{1}}+\lambda \nu \psi ^{\nu
-1}u\right) ^{2}+\dsum\limits_{i,j=2}^{n}u_{x_{i}}^{2}\right) \varphi
_{\lambda ,\nu }\psi ^{-\nu +1}\right] + \\ 
+\partial _{t}\left( -\lambda ^{2}\nu ^{2}\psi ^{\nu -1}\left( 1-2\psi
^{-\nu }\left( \nu -1\right) /\left( \lambda \nu \right) \right)
u^{2}\varphi _{\lambda ,\nu }+\partial _{t}\left( \left( \lambda /2\right)
u^{2}\varphi _{\lambda ,\nu }\right) \right) ,%
\end{array}%
\right.  \label{7.26}
\end{equation}%
\begin{equation}
\left. 
\begin{array}{c}
\func{div}U_{6}=\func{div}U_{4}+\func{div}\left( \lambda U_{5}\right) = \\ 
\begin{array}{c}
=\left( -2u_{t}\left( u_{x_{1}}+\lambda \nu \psi ^{\nu -1}u\right) \varphi
_{\lambda ,\nu }\psi ^{-\nu +1}\right) _{x_{1}} \\ 
+\left( -2\lambda \nu \left( u_{x_{1}}+\lambda \nu \psi ^{\nu -1}u\right)
^{2}\varphi _{\lambda ,\nu }+2\lambda \nu
\dsum\limits_{i=2}^{n}u_{x_{i}}^{2}\varphi _{\lambda ,\nu }\right) _{x_{1}}+
\\ 
+\left( -2\lambda ^{3}\nu ^{3}\psi ^{2\nu -2}\left( 1-2\psi ^{-\nu }\left(
\nu -1\right) /\left( \lambda \nu \right) \right) u^{2}\varphi _{\lambda
,\nu }\right) _{x_{1}} \\ 
+\dsum\limits_{i=2}^{n}\left( -4\lambda \nu \left( u_{x_{1}}+\lambda \nu
\psi ^{\nu -1}u\right) u_{x_{i}}\varphi _{\lambda ,\nu
}-2u_{t}u_{x_{i}}\varphi _{\lambda ,\nu }\psi ^{-\nu +1}\right) _{x_{i}}+ \\ 
+\left( -\lambda u_{x_{1}}u\varphi _{\lambda ,\nu }+\lambda ^{2}\nu \psi
^{\nu -1}u^{2}\varphi _{\lambda ,\nu }\right)
_{x_{1}}+\dsum\limits_{i=2}^{n}\left( -\lambda u_{x_{i}}u\varphi _{\lambda
,\nu }\right) _{x_{i}}.%
\end{array}%
\end{array}%
\right.  \label{7.27}
\end{equation}%
Applying Cauchy-Schwarz inequality to the left hand side of (\ref{7.25}) and
also using (\ref{7.1}), we obtain a pointwise Carleman estimate from the
above for the lower order derivatives via $\left( u_{t}-\Delta u\right)
^{2}\varphi _{\lambda ,\nu },$%
\begin{equation}
\left. 
\begin{array}{c}
\left( u_{t}-\Delta u\right) ^{2}\varphi _{\lambda ,\nu }\geq C\lambda
\left( \nabla u\right) ^{2}\varphi _{\lambda ,\nu }+C\lambda ^{3}\nu
^{4}\psi ^{2\nu -2}u^{2}+\partial _{t}V_{5}+\func{div}U_{6}, \\ 
\forall \lambda \geq \lambda _{0}=\lambda _{0}\left( A_{1}\right) >1,\text{ }%
\nu =\nu _{0}=\nu _{0}\left( A_{1}\right) >2.%
\end{array}%
\right.  \label{7.28}
\end{equation}

We now need to incorporate estimates for derivatives $u_{x_{i}x_{j}}$ and $%
u_{t}$ in a close analog of (\ref{7.28}).

\subsubsection{Step 8. Estimate again $\left( u_{t}-\Delta u\right) ^{2}%
\protect\varphi _{\protect\lambda ,\protect\nu _{0}}$ from the below}

We now set $\nu =\nu _{0},$see the beginning of the proof of Theorem 4.1. We
have%
\begin{equation}
\left( u_{t}-\Delta u\right) ^{2}\varphi _{\lambda ,\nu
_{0}}=u_{t}^{2}\varphi _{\lambda ,\nu _{0}}-2u_{t}\Delta u\varphi _{\lambda
,\nu _{0}}+\left( \Delta u\right) ^{2}\varphi _{\lambda ,\nu _{0}}.
\label{7.29}
\end{equation}%
We estimate separately from the below the second and the third terms in the
righ hand side of (\ref{7.29}).

\textbf{Step 8.1.} First, estimate from the below the term $u_{t}^{2}\varphi
_{\lambda ,\nu _{0}}-2u_{t}\Delta u\varphi _{\lambda ,\nu _{0}}$ in (\ref%
{7.29}).

We have:%
\begin{equation}
\left. 
\begin{array}{c}
-2u_{t}\Delta u\varphi _{\lambda ,\nu _{0}}=-2u_{t}u_{x_{1}x_{1}}\varphi
_{\lambda ,\nu _{0}}+\dsum\limits_{i=2}^{n}\left(
-2u_{t}u_{x_{i}x_{i}}\varphi _{\lambda ,\nu _{0}}\right) = \\ 
=\left( -2u_{t}u_{x_{1}}\varphi _{\lambda ,\nu _{0}}\right)
_{x_{1}}+2u_{x_{1}t}u_{x_{1}}\varphi _{\lambda ,\nu _{0}}+2\lambda \nu
_{0}\psi ^{\nu _{0}-1}u_{t}u_{x_{1}}\varphi _{\lambda ,\nu _{0}}+ \\ 
+\dsum\limits_{i=2}^{n}\left( -2u_{t}u_{x_{i}}\varphi _{\lambda ,\nu
_{0}}\right) _{x_{i}}+\dsum\limits_{i=2}^{n}2u_{x_{i}t}u_{x_{i}}\varphi
_{\lambda ,\nu _{0}}= \\ 
=\partial _{t}\left( \left( \nabla u\right) ^{2}\varphi _{\lambda ,\nu
_{0}}\right) +\dsum\limits_{i=1}^{n}\left( -2u_{t}u_{x_{i}}\varphi _{\lambda
,\nu _{0}}\right) _{x_{i}}+2\lambda \nu _{0}\psi ^{\nu
_{0}-1}u_{t}u_{x_{1}}\varphi _{\lambda ,\nu _{0}}.%
\end{array}%
\right.  \label{7.30}
\end{equation}%
By Young's inequality and (\ref{4.2}) 
\begin{equation}
\left. 
\begin{array}{c}
2\lambda \nu _{0}\psi ^{\nu _{0}-1}u_{t}u_{x_{1}}\varphi _{\lambda ,\nu
_{0}}\geq -\left( 1/2\right) u_{t}^{2}\varphi _{\lambda ,\nu _{0}}-2\lambda
^{2}\nu _{0}^{2}\psi ^{2\nu _{0}-2}u_{x_{1}}^{2}\varphi _{\lambda ,\nu
_{0}}\geq \\ 
\geq -\left( 1/2\right) u_{t}^{2}\varphi _{\lambda ,\nu _{0}}-C\lambda
^{2}u_{x_{1}}^{2}\varphi _{\lambda ,\nu _{0}}.%
\end{array}%
\right.  \label{7.31}
\end{equation}%
Hence, using (\ref{7.29})-(\ref{7.31}), we obtain%
\begin{equation}
\left. 
\begin{array}{c}
\left( u_{t}-\Delta u\right) ^{2}\varphi _{\lambda ,\nu _{0}}\geq \left(
1/2\right) u_{t}^{2}\varphi _{\lambda ,\nu _{0}}-C\lambda
^{2}u_{x_{1}}^{2}\varphi _{\lambda ,\nu _{0}}+\left( \Delta u\right)
^{2}\varphi _{\lambda ,\nu _{0}}+ \\ 
+\func{div}U_{7},%
\end{array}%
\right.  \label{7.32}
\end{equation}%
\begin{equation}
\partial _{t}V_{6}=\partial _{t}\left( \left( \nabla u\right) ^{2}\varphi
_{\lambda ,\nu _{0}}\right) ,  \label{7.33}
\end{equation}%
\begin{equation}
\func{div}U_{7}=\dsum\limits_{i=1}^{n}\left( -2u_{t}u_{x_{i}}\varphi
_{\lambda ,\nu _{0}}\right) _{x_{i}}.  \label{7.34}
\end{equation}

\textbf{Step 8.2.} Second, estimate from the term $\left( \Delta u\right)
^{2}\varphi _{\lambda ,\nu _{0}}$ in (\ref{7.29}).

We have:%
\begin{equation}
\left( \Delta u\right) ^{2}\varphi _{\lambda ,\nu
_{0}}=u_{x_{1}x_{1}}^{2}\varphi _{\lambda ,\nu
_{0}}+\dsum\limits_{i=2}^{n}2u_{x_{1}x_{1}}u_{x_{i}x_{i}}\varphi _{\lambda
,\nu _{0}}+\dsum\limits_{i,j=2}^{n}u_{x_{i}x_{i}}u_{x_{j}x_{j}}\varphi
_{\lambda ,\nu _{0}}.  \label{7.35}
\end{equation}

Estimate the term:%
\begin{equation}
\left. 
\begin{array}{c}
\dsum\limits_{i=2}^{n}2u_{x_{1}x_{1}}u_{x_{i}x_{i}}\varphi _{\lambda ,\nu
_{0}}=\dsum\limits_{i=2}^{n}\left( 2u_{x_{1}x_{1}}u_{x_{i}}\varphi _{\lambda
,\nu _{0}}\right) _{x_{i}}-\dsum\limits_{i=2}^{n}\left(
2u_{x_{1}x_{1}x_{i}}u_{x_{i}}\varphi _{\lambda ,\nu _{0}}\right) = \\ 
=\dsum\limits_{i=2}^{n}\left( 2u_{x_{1}x_{1}}u_{x_{i}}\varphi _{\lambda ,\nu
_{0}}\right) _{x_{i}}+\dsum\limits_{i=2}^{n}\left(
-2u_{x_{1}x_{i}}u_{x_{i}}\varphi _{\lambda ,\nu _{0}}\right)
_{x_{1}}+2\dsum\limits_{i=2}^{n}u_{x_{1}x_{i}}^{2}\varphi _{\lambda ,\nu
_{0}}+ \\ 
+\dsum\limits_{i=2}^{n}4\lambda \nu _{0}\psi ^{\nu
_{0}-1}u_{x_{1}x_{i}}u_{x_{i}}\varphi _{\lambda ,\nu _{0}}.%
\end{array}%
\right.  \label{7.36}
\end{equation}%
By Young's inequality and (\ref{4.2})%
\begin{equation*}
\dsum\limits_{i=2}^{n}4\lambda \nu _{0}\psi ^{\nu
_{0}-1}u_{x_{1}x_{i}}u_{x_{i}}\varphi _{\lambda ,\nu _{0}}\geq
-\dsum\limits_{i=2}^{n}u_{x_{1}x_{i}}^{2}\varphi _{\lambda ,\nu
_{0}}-C\lambda ^{2}\left( \nabla u\right) ^{2}\varphi _{\lambda ,\nu _{0}}.
\end{equation*}%
Hence, using (\ref{7.35}) and (\ref{7.36}), we obtain%
\begin{equation}
\left( \Delta u\right) ^{2}\varphi _{\lambda ,\nu _{0}}\geq
\dsum\limits_{i=1}^{n}u_{x_{1}x_{i}}^{2}\varphi _{\lambda ,\nu
_{0}}+\dsum\limits_{i,j=2}^{n}u_{x_{i}x_{i}}u_{x_{j}x_{j}}\varphi _{\lambda
,\nu _{0}}-C\lambda ^{2}\left( \nabla u\right) ^{2}\varphi _{\lambda ,\nu
_{0}}+\func{div}U_{8},  \label{7.37}
\end{equation}%
\begin{equation}
\func{div}U_{8}=\dsum\limits_{i=2}^{n}\left(
-2u_{x_{1}x_{i}}u_{x_{i}}\varphi _{\lambda ,\nu _{0}}\right)
_{x_{1}}+\dsum\limits_{i=2}^{n}\left( 2u_{x_{1}x_{1}}u_{x_{i}}\varphi
_{\lambda ,\nu _{0}}\right) _{x_{i}}.  \label{7.38}
\end{equation}

We now estimate the following term in (\ref{7.37}):%
\begin{equation*}
\left. 
\begin{array}{c}
\dsum\limits_{i,j=2}^{n}u_{x_{i}x_{i}}u_{x_{j}x_{j}}\varphi _{\lambda ,\nu
_{0}}=\dsum\limits_{i,j=2}^{n}\left( u_{x_{i}x_{i}}u_{x_{j}}\varphi
_{\lambda ,\nu _{0}}\right)
_{x_{j}}-\dsum\limits_{i,j=2}^{n}u_{x_{i}x_{i}x_{j}}u_{x_{j}}\varphi
_{\lambda ,\nu _{0}}= \\ 
=\dsum\limits_{i,j=2}^{n}\left( u_{x_{i}x_{i}}u_{x_{j}}\varphi _{\lambda
,\nu _{0}}\right) _{x_{j}}+\dsum\limits_{i,j=2}^{n}\left(
-u_{x_{i}x_{j}}u_{x_{j}}\varphi _{\lambda ,\nu _{0}}\right)
_{x_{i}}+\dsum\limits_{i,j=2}^{n}u_{x_{i}x_{j}}^{2}\varphi _{\lambda ,\nu
_{0}}.%
\end{array}%
\right.
\end{equation*}%
Combining this with (\ref{7.37}) and (\ref{7.38}), we obtain%
\begin{equation}
\left( \Delta u\right) ^{2}\varphi _{\lambda ,\nu _{0}}\geq
\dsum\limits_{i=1}^{n}u_{x_{i}x_{j}}^{2}\varphi _{\lambda ,\nu
_{0}}-C\lambda ^{2}\left( \nabla u\right) ^{2}\varphi _{\lambda ,\nu _{0}}+%
\func{div}U_{9},  \label{7.39}
\end{equation}%
\begin{equation}
\left. 
\begin{array}{c}
\func{div}U_{9}=\dsum\limits_{i=2}^{n}\left(
-2u_{x_{1}x_{i}}u_{x_{i}}\varphi _{\lambda ,\nu _{0}}\right)
_{x_{1}}+\dsum\limits_{i=2}^{n}\left( 2u_{x_{1}x_{1}}u_{x_{i}}\varphi
_{\lambda ,\nu _{0}}\right) _{x_{i}}+ \\ 
+\dsum\limits_{i,j=2}^{n}\left( u_{x_{i}x_{i}}u_{x_{j}}\varphi _{\lambda
,\nu _{0}}\right) _{x_{j}}+\dsum\limits_{i,j=2}^{n}\left(
-u_{x_{i}x_{j}}u_{x_{j}}\varphi _{\lambda ,\nu _{0}}\right) _{x_{i}}.%
\end{array}%
\right.  \label{7.40}
\end{equation}%
Combining (\ref{7.39}) and (\ref{7.40}) with (\ref{7.32})-(\ref{7.34}), we
obtain%
\begin{equation}
\left. 
\begin{array}{c}
\left( u_{t}-\Delta u\right) ^{2}\varphi _{\lambda ,\nu _{0}}\geq \left(
1/2\right) u_{t}^{2}\varphi _{\lambda ,\nu
_{0}}+\dsum\limits_{i=1}^{n}u_{x_{i}x_{j}}^{2}\varphi _{\lambda ,\nu
_{0}}-C\lambda ^{2}\left( \nabla u\right) ^{2}\varphi _{\lambda ,\nu _{0}}+
\\ 
+\partial _{t}\left( \left( \nabla u\right) ^{2}\varphi _{\lambda ,\nu
_{0}}\right) +\func{div}U_{10},%
\end{array}%
\right.  \label{7.41}
\end{equation}%
\begin{equation}
\left. 
\begin{array}{c}
\func{div}U_{10}=\func{div}U_{7}+\func{div}U_{9}= \\ 
=\dsum\limits_{i=1}^{n}\left( -2u_{t}u_{x_{i}}\varphi _{\lambda ,\nu
_{0}}\right) _{x_{i}}+ \\ 
\begin{array}{c}
+\dsum\limits_{i=2}^{n}\left( -2u_{x_{1}x_{i}}u_{x_{i}}\varphi _{\lambda
,\nu _{0}}\right) _{x_{1}}+\dsum\limits_{i=2}^{n}\left(
2u_{x_{1}x_{1}}u_{x_{i}}\varphi _{\lambda ,\nu _{0}}\right) _{x_{i}}+ \\ 
+\dsum\limits_{i,j=2}^{n}\left( u_{x_{i}x_{i}}u_{x_{j}}\varphi _{\lambda
,\nu _{0}}\right) _{x_{j}}+\dsum\limits_{i,j=2}^{n}\left(
-u_{x_{i}x_{j}}u_{x_{j}}\varphi _{\lambda ,\nu _{0}}\right) _{x_{i}}.%
\end{array}%
\end{array}%
\right.  \label{7.42}
\end{equation}

\subsubsection{Step 9. Divide (\protect\ref{7.41}) and (\protect\ref{7.42})
by $2\protect\lambda $ and sum up with (\protect\ref{7.28}), taking into
account (\protect\ref{7.26}) and (\protect\ref{7.27})}

After suming up as indicated, we divide both parts of the resulting estimate
by $\left( 1+1/\left( 2\lambda \right) \right) .$ Since we use the constant $%
C$, then this division will affect only terms under $\partial _{t}$ and $%
\func{div}$ signs. We obtain the target estimate (\ref{4.3}) of this
theorem, where%
\begin{equation}
\partial _{t}V=\partial _{t}\left( \frac{2\lambda V_{5}+\left( \nabla
u\right) ^{2}\varphi _{\lambda ,\nu _{0}}}{2\lambda +1}\right) ,\func{div}U=%
\func{div}\left( \frac{2\lambda U_{6}+U_{10}}{2\lambda +1}\right) .
\label{7.43}
\end{equation}%
Formulas (\ref{7.43}) are equivalent with formulas (\ref{4.30}), (\ref{4.31}%
). \ $\square $

\subsection{Proof of Theorem 4.2}

\label{sec:6.2}

It is convenient to assume first that $u\in C^{4,2}\left( \overline{Q}%
_{T}\right) $ since Theorem 4.1 is proven for these functions. Integrate (%
\ref{4.3}) over $Q_{T}.$ It follows from (\ref{4.30}), (\ref{4.5}) and the
last line of (\ref{4.6}) that integrals over $\left\{ t=0\right\} $ and $%
\left\{ t=T\right\} $ are mutually canceled. Therefore, by Gauss formula 
\begin{equation}
\left. 
\begin{array}{c}
\dint\limits_{Q_{T}}\left( u_{t}-\Delta u\right) ^{2}\varphi _{\lambda ,\nu
_{0}}dxdt\geq \left( C/\lambda \right) \dint\limits_{Q_{T}}\left(
u_{t}^{2}+\dsum\limits_{i,j=1}^{n}u_{x_{i}x_{j}}^{2}\right) \varphi
_{\lambda ,\nu _{0}}dxdt+ \\ 
+C\dint\limits_{Q_{T}}\left[ \lambda \left( \nabla u\right) ^{2}+\lambda
^{3}u^{2}\right] \varphi _{\lambda ,\nu _{0}}dxdt+\dint\limits_{S_{T}}U\cos
\left( \mu ,x\right) dS,\text{ }\forall \lambda \geq \lambda _{0},%
\end{array}%
\right.  \label{7.44}
\end{equation}%
where $\mu $ is the outward looking unit normal vector at $\partial \Omega $.

We now evaluate the term%
\begin{equation}
\dint\limits_{S_{T}}U\cos \left( \mu ,x\right) dS.  \label{7.45}
\end{equation}%
To do this, we use (\ref{2.101})-(\ref{2.105}) and (\ref{4.31}). First,
consider the part $\Gamma _{T}^{+}$ of $S_{T}.$ Obviously, $\mu =\left(
1,0,...,0\right) $ on $\Gamma _{T}^{+}.$ Note that $dS=dx_{2}...dx_{n}dt$ on 
$\Gamma _{T}^{+}.$ By (\ref{4.31}) and (\ref{7.45}) 
\begin{equation*}
\left. 
\begin{array}{c}
\dint\limits_{\Gamma _{T}^{+}}U\cos \left( \mu ,x\right) dS= \\ 
=\dint\limits_{\Gamma _{T}^{+}}\left[ \left( 2\lambda /\left( 2\lambda
+1\right) \right) \left( -2u_{t}\left( u_{x_{1}}+\lambda \nu _{0}\psi ^{\nu
_{0}-1}u\right) \varphi _{\lambda ,\nu _{0}}\psi ^{-\nu _{0}+1}\right) %
\right] dS+ \\ 
+\dint\limits_{\Gamma _{T}^{+}}\left[ \left( 2\lambda /\left( 2\lambda
+1\right) \right) \left( -2\lambda \nu _{0}\left( u_{x_{1}}+\lambda \nu
_{0}\psi ^{\nu _{0}-1}u\right) ^{2}\varphi _{\lambda ,\nu _{0}}\right) %
\right] dS+ \\ 
+\dint\limits_{\Gamma _{T}^{+}}\left[ \left( 2\lambda /\left( 2\lambda
+1\right) \right) \left( 2\lambda \nu
_{0}\dsum\limits_{i=2}^{n}u_{x_{i}}^{2}\varphi _{\lambda ,\nu _{0}}\right) %
\right] dS+ \\ 
+\dint\limits_{\Gamma _{T}^{+}}\left[ \left( 2\lambda /\left( 2\lambda
+1\right) \right) \left( -2\lambda ^{3}\nu _{0}^{3}\psi ^{2\nu _{0}-2}\left(
1-2\psi ^{-\nu _{0}}\left( \nu _{0}-1\right) /\left( \lambda \nu _{0}\right)
\right) u^{2}\varphi _{\lambda ,\nu _{0}}\right) \right] dS+ \\ 
+\dint\limits_{\Gamma _{T}^{+}}\left[ \left( 2\lambda /\left( 2\lambda
+1\right) \right) \left( -\lambda u_{x_{1}}u\varphi _{\lambda ,\nu
_{0}}+\lambda ^{2}\nu _{0}\psi ^{\nu _{0}-1}u^{2}\varphi _{\lambda ,\nu
_{0}}\right) \right] dS+ \\ 
+\dint\limits_{\Gamma _{T}^{+}}\dsum\limits_{i=2}^{n}\left[ \left( 1/\left(
2\lambda +1\right) \right) \left( -2u_{x_{1}x_{i}}u_{x_{i}}\varphi _{\lambda
,\nu _{0}}\right) \right] dS.%
\end{array}%
\right.
\end{equation*}%
Combining this equality with (\ref{4.2}), we obtain%
\begin{equation}
\left. 
\begin{array}{c}
\dint\limits_{\Gamma _{T}^{+}}U\cos \left( n,x\right) dS\geq \\ 
\geq -C\exp \left[ 3\lambda \left( 2A_{1}+2\right) ^{\nu _{0}}\right] \left(
\left\Vert u\right\Vert _{H^{2,1}\left( \Gamma _{T}^{+}\right)
}^{2}+\left\Vert u_{x_{1}}\right\Vert _{H^{1,0}\left( \Gamma _{T}^{+}\right)
}^{2}\right) .%
\end{array}%
\right.  \label{7.46}
\end{equation}%
Second, consider the part $\Gamma _{T}^{+}$ of $S_{T}.$ by (\ref{2.1}) and (%
\ref{4.1}) $\varphi _{\lambda ,\nu _{0}}=\exp \left( 2^{\nu _{0}+1}\lambda
\right) $ on $\Gamma _{T}^{-}.$ Hence, we obtain similarly with (\ref{7.46})%
\begin{equation}
\left. 
\begin{array}{c}
\dint\limits_{\Gamma _{T}^{-}}U\cos \left( \mu ,x\right) dS\geq \\ 
\geq -C\exp \left[ 3\cdot 2^{\nu _{0}}\lambda \right] \left( \left\Vert
u\right\Vert _{H^{2,1}\left( \Gamma _{T}^{-}\right) }^{2}+\left\Vert
u_{x_{1}}\right\Vert _{H^{1,0}\left( \Gamma _{T}^{-}\right) }^{2}\right) .%
\end{array}%
\right.  \label{7.47}
\end{equation}

We now evaluate the integral%
\begin{equation}
\dint\limits_{\partial _{i}^{+}\Omega _{T}}U\cos \left( \mu ,x\right) dS,%
\text{ }i=2,...,n.  \label{7.48}
\end{equation}%
Obviously $\mu =\left( \delta _{i1},...,\delta _{in}\right) ,$ where 
\begin{equation*}
\delta _{ij}=\left\{ 
\begin{array}{c}
1\text{ if }i=j, \\ 
0\text{ if }i\neq j.%
\end{array}%
\right.
\end{equation*}%
Hence, by (\ref{4.31}) and (\ref{7.48}) for $i\in \left[ 2,n\right] $%
\begin{equation}
\left. 
\begin{array}{c}
\dint\limits_{\partial _{i}^{+}\Omega _{T}}U\cos \left( \mu ,x\right) dS= \\ 
+\dint\limits_{\partial _{i}^{+}\Omega _{T}}\left[ \left( 2\lambda /\left(
2\lambda +1\right) \right) \left( -4\lambda \nu _{0}\left( u_{x_{1}}+\lambda
\nu _{0}\psi ^{\nu _{0}-1}u\right) u_{x_{i}}\varphi _{\lambda ,\nu
_{0}}\right) \right] dS+ \\ 
+\dint\limits_{\partial _{i}^{+}\Omega _{T}}\left[ \left( 2\lambda /\left(
2\lambda +1\right) \right) \left( -2u_{t}u_{x_{i}}\right) \psi ^{-\nu
_{0}+1}\varphi _{\lambda ,\nu _{0}}\right] dS+ \\ 
+\dint\limits_{\partial _{i}^{+}\Omega _{T}}\left[ \left( 2\lambda /\left(
2\lambda +1\right) \right) \left( -\lambda u_{x_{i}}u-2u_{t}u_{x_{i}}\right)
\varphi _{\lambda ,\nu _{0}}\right] dS+ \\ 
\dint\limits_{\partial _{i}^{+}\Omega _{T}}\left[ \left( 1/\left( 2\lambda
+1\right) \right) \left( u_{x_{j}x_{j}}u_{x_{i}}\varphi _{\lambda ,\nu
_{0}}-u_{x_{i}x_{j}}u_{x_{j}}\varphi _{\lambda ,\nu _{0}}\right) \right] dS.%
\end{array}%
\right.  \label{7.49}
\end{equation}%
Suppose that $i=j$ in the last line of (\ref{7.49}). Then this term is:%
\begin{equation}
\left. 
\begin{array}{c}
\dint\limits_{\partial _{i}^{+}\Omega _{T}}\left[ \left( 1/\left( 2\lambda
+1\right) \right) \left( u_{x_{j}x_{j}}u_{x_{i}}\varphi _{\lambda ,\nu
_{0}}-u_{x_{i}x_{j}}u_{x_{j}}\varphi _{\lambda ,\nu _{0}}\right) \right] dS=
\\ 
=\dint\limits_{\partial _{i}^{+}\Omega _{T}}\left[ \left( 1/\left( 2\lambda
+1\right) \right) \left( u_{x_{i}x_{i}}u_{x_{i}}\varphi _{\lambda ,\nu
_{0}}-u_{x_{i}x_{i}}u_{x_{i}}\varphi _{\lambda ,\nu _{0}}\right) \right]
dS=0.%
\end{array}%
\right.  \label{7.50}
\end{equation}%
Hence, the last line of (\ref{7.49}) is not identically zero only if $i\neq
j $. Hence, by (\ref{4.2})%
\begin{equation}
\left. 
\begin{array}{c}
\dint\limits_{\partial _{i}^{+}\Omega _{T}}\left[ \left( 1/\left( 2\lambda
+1\right) \right) \left( u_{x_{j}x_{j}}u_{x_{i}}\varphi _{\lambda ,\nu
_{0}}-u_{x_{i}x_{j}}u_{x_{j}}\varphi _{\lambda ,\nu _{0}}\right) \right]
dS\geq \\ 
\geq -C\exp \left[ 3\lambda \left( 2A_{1}+2\right) ^{\nu _{0}}\right] \left(
\left\Vert u\right\Vert _{H^{2,1}\left( \partial _{i}^{+}\Omega _{T}\right)
}^{2}+\left\Vert \partial _{n}u\right\Vert _{H^{1,0}\left( \partial
_{i}^{+}\Omega _{T}\right) }^{2}\right) .%
\end{array}%
\right.  \label{7.51}
\end{equation}%
Analysis of the sum of second, third and fourth lines of (\ref{7.49}) shows
that this sum can also be estimated from the below like in (\ref{7.51}).
Thus, combining the latter considerations with (\ref{7.45})-(\ref{7.51}), we
obtain%
\begin{equation*}
\dint\limits_{S_{T}}U\cos \left( \mu ,x\right) dS\geq -C\exp \left[ 3\lambda
\left( 2A_{1}+2\right) ^{\nu _{0}}\right] \left( \left\Vert u\right\Vert
_{H^{2,1}\left( S_{T}\right) }^{2}+\left\Vert \partial _{n}u\right\Vert
_{H^{1,0}\left( S_{T}\right) }^{2}\right) ,
\end{equation*}%
and also, more precisely, by (\ref{7.47}),%
\begin{equation}
\left. 
\begin{array}{c}
\dint\limits_{S_{T}}U\cos \left( \mu ,x\right) dS\geq \\ 
\geq -C\exp \left[ 3\lambda \left( 2A_{1}+2\right) ^{\nu _{0}}\right] \left(
\left\Vert u\right\Vert _{H^{2,1}\left( S_{T}\diagdown \Gamma
_{T}^{-}\right) }^{2}+\left\Vert \partial _{n}u\right\Vert _{H^{1,0}\left(
\Gamma _{T}^{-}\right) }^{2}\right) - \\ 
-C\exp \left[ 3\cdot 2^{\nu _{0}}\lambda \right] \left( \left\Vert
u\right\Vert _{H^{2,1}\left( \Gamma _{T}^{-}\right) }^{2}+\left\Vert
u_{x_{1}}\right\Vert _{H^{1,0}\left( \Gamma _{T}^{-}\right) }^{2}\right) .%
\text{ }%
\end{array}%
\right.  \label{7.54}
\end{equation}%
Finally, density arguments ensure that one can replace $u\in C^{4,2}\left( 
\overline{Q}_{T}\right) $ with $u\in H^{4,2}\left( Q_{T}\right) $ in (\ref%
{7.44}) and (\ref{7.54}). The target estimate (\ref{4.6}) of this theorem
follows immediately. \ $\square $

\end{document}